\documentclass[aps,pre,10pt,showpacs,floatfix,onecolumn,nofootinbib,a4paper]{revtex4-2}

\newtheorem{theorem}{\bf Theorem}[section]
\newtheorem{condition}{\bf Condition}[section]

\usepackage{amssymb, amsmath, amsthm, mathtools}
\usepackage{bm}
\usepackage{mathrsfs}
\usepackage{hyperref,bookmark}
\usepackage{graphicx}
\usepackage{epsfig,color}
\usepackage{float}
\usepackage{subcaption}
\usepackage{verbatim}
\usepackage{dcolumn}
\usepackage{tikz}
\usepackage{url}

\usepackage[utf8x]{inputenc}
\usepackage[english]{babel}
\usepackage{mathrsfs}
\usepackage{epsfig}
\usepackage{listings}

\usepackage{mdwlist}
\usepackage{bm}
\usepackage{dsfont}
\usepackage{oldgerm}
\usepackage{geometry}
\usepackage{relsize}
\usepackage{afterpage}
\usepackage{lmodern}
\usepackage{tabularx}
\usepackage{cancel}
\usepackage{dutchcal} 

\newcommand{\dd}{\mathrm{d}}
\newcommand{\ex}{\mathrm{e}}

\newcommand{\vt}{\vartheta}
\newcommand{\vp}{\varphi}
\newcommand{\trans}{^\mathsf{T}}

\newcommand{\av}{\bm{a}}
\newcommand{\n}{\bm{n}}

\newcommand{\cv}{\bm{c}}

\newcommand{\dv}{\bm{d}}
\newcommand{\df}{\dv_{\f}}

\newcommand{\e}{\bm{e}}
\newcommand{\x}{\bm{x}}
\newcommand{\zero}{\bm{0}}
\newcommand{\y}{\bm{y}}

\newcommand{\normal}{\bm{\nu}}
\newcommand{\surface}{\mathscr{S}}
\newcommand{\tplane}{\mathscr{T}}

\newcommand{\nablas}{\nabla\!_\mathrm{s}}
\newcommand{\lap}{\mathop{{}\bigtriangleup}\nolimits}
\newcommand{\laps}{\mathop{{}\bigtriangleup_\mathrm{s}\!}\nolimits}
\newcommand{\nablastwo}{\nabla^2\!\!\!\!_\mathrm{s}}
\newcommand{\nablan}{\nabla\!_\nu}

\newcommand{\uv}{\bm{u}}

\newcommand{\f}{\bm{f}}
\newcommand{\g}{\bm{g}}
\newcommand{\h}{\bm{h}}

\newcommand{\rv}{\bm{r}}

\newcommand{\bv}{\bm{b}}
\newcommand{\avf}{\av_{\f}}
\newcommand{\bvf}{\bv_{\f}}
\newcommand{\bg}{\bv_{\g}}

\newcommand{\vv}{\bm{v}}
\newcommand{\con}{\bm{c}}
\newcommand{\W}{\mathbf{W}}
\newcommand{\F}{\mathbf{F}}

\newcommand{\I}{\mathbf{I}}

\newcommand{\R}{\mathbf{R}}
\newcommand{\Rf}{\R_{\f}}
\newcommand{\Rg}{\R_{\g}}
\newcommand{\Ry}{\R_{\y}}
\newcommand{\Rfb}{\Rf^{(\mathrm{b})}}

\newcommand{\Ryb}{\Ry^{(\mathrm{b})}}

\newcommand{\Rfd}{\Rf^{(\mathrm{d})}}

\newcommand{\C}{\mathbf{C}}
\newcommand{\B}{\mathbf{B}}
\newcommand{\U}{\mathbf{U}}
\newcommand{\V}{\mathbf{V}}
\newcommand{\Cf}{\C_{\f}}
\newcommand{\Bf}{\B_{\f}}
\newcommand{\Uf}{\U_{\f}}
\newcommand{\Vf}{\V_{\f}}

\newcommand{\Uy}{\U_{\y}}

\newcommand{\Lt}{\mathbf{L}}
\newcommand{\M}{\mathbf{M}}
\newcommand{\curl}{\operatorname{curl}}
\newcommand{\curls}{\curl_\mathrm{s}\!}
\newcommand{\divs}{\diver_\mathrm{s}}
\newcommand{\tr}{\operatorname{tr}}
\newcommand{\co}{\operatorname{co}}
\newcommand{\skw}{\operatorname{skw}}

\newcommand{\nay}{(\nabla\y)}
\newcommand{\naf}{(\nabla\f)}

\newcommand{\curvature}{(\nablas\normal)}
\newcommand{\curvaturast}{(\nablas\normal^\ast)}

\newcommand{\euclid}{\mathscr{E}}
\newcommand{\translation}{\mathscr{V}}

\newcommand{\framen}{(\n_1,\n_2,\normal)}
\newcommand{\framee}{(\e_1,\e_2,\e_3)}

\newcommand{\tangent}{\bm{t}}
\newcommand{\curve}{\mathcal{c}}
\newcommand{\subsurface}{\mathcal{s}}

\newcommand{\Rnu}{\mathbf{R}_{\normal}}
\newcommand{\pure}{\mathbf{A}}
\newcommand{\orth}{\mathsf{SO}(3)}
\newcommand{\orthnu}{\mathsf{SO}(\normal)}
\newcommand{\proj}{\mathbf{P}(\normal)}
\newcommand{\proje}{\mathbf{P}(\e)}
\newcommand{\Proj}{\mathbf{P}(\e_3)}
\newcommand{\ball}{\mathbb{B}}
\newcommand{\sphere}{\mathbb{S}^2}

\renewcommand{\leq}{\leqslant}

\newcommand{\diver}{\operatorname{div}} 

\newcommand{\nigh}[1]{{\color{black}{#1}}}

\begin{document}
\title{Bending-Neutral Deformations of Minimal Surfaces}
\author{Andr\'e M. Sonnet}
\email{Andre.Sonnet@strath.ac.uk}
\affiliation{Department of Mathematics and Statistics, University of Strathclyde, 26 Richmond Street, Glasgow, G1 1XH, U.K. }
\author{Epifanio G. Virga}
\email{eg.virga@unipv.it}
\affiliation{Dipartimento di Matematica, Universit\`a di Pavia, Via Ferrata 5, 27100 Pavia, Italy }

\date{\today}

\begin{abstract}
	Minimal surfaces are ubiquitous in nature. Here they are considered as geometric objects that bear a deformation content. By refining the resolution of the surface deformation gradient afforded by the polar decomposition theorem, we identify a \emph{bending} content and a class of deformations that leave it unchanged. These are the \emph{bending-neutral} deformations, fully characterized by an integrability condition; \nigh{they preserve normals}. We prove that (1) every minimal surface is transformed into a minimal surface by a bending-neutral deformation, (2) given two minimal surfaces \nigh{with the same system of normals}, there is a bending-neutral deformation that maps one into the other, \nigh{and (3)} \emph{all}  minimal surfaces have indeed a \emph{universal} bending content.
\end{abstract}

\maketitle
\section{Introduction}\label{sec:intro}
The \emph{universe is parsimonious}: so says the title of a celebrated book \cite{hildebrandt:parsimonious}, which among many other things contains a fascinating account on minimal surfaces (disguised as soap bubbles). The theory of plates and shells should make no exception to this universal rule. However, researchers are still debating the most appropriate kinematic measures of deformation that would enter the elastic energy functional in an \emph{intrinsic} (direct) theory of these bodies. In the language of the book \cite{antman:nonlinear}, an intrinsic theory represents plates and shells as truly two-dimensional bodies, with mass distributed on a surface $\surface$ in three-dimensional Euclidean space $\euclid$. The following works provide but a partial sample of the current debate \cite{ghiba:isotropic,ghiba:constrained,vitral:dilation,vitral:energies,virga:pure,acharya:mid-surface,ghiba:essay,vitral:assorted}.

Here we build on an earlier proposal for measures of \emph{pure} bending of a surface $\surface$ \cite{virga:pure}, which was formulated as an invariance requirement under the class of \emph{bending-neutral} deformations, \nigh{which preserve normals}. In this paper, we identify as \emph{bending content} of $\surface$ the geometric object (a vector field) left unchanged by these deformations.

We prove that \emph{all} minimal surfaces  share \emph{one} and the same bending content \nigh{and those that share the same system of normals}  can be seen as the image of one another under a bending-neutral deformation. Thus, making it impossible to properly \emph{bend} one minimal surface into another, \nigh{if the system of normals is to be \emph{preserved}.}

The paper is organized as follows. In Sect.~\ref{sec:calculus}, to make our development self-contained, we collect a few results about surface tensor calculus that are used in the rest of the paper. In Sect.~\ref{sec:deformations}, we lay out our kinematic analysis: we recall the definition of bending-neutral deformations and identify the bending content of a surface. In Sect.~\ref{sec:minimal}, we show how minimal surfaces are tightly connected by bending-neutral deformations; \nigh{we accordingly introduce the notion of bending-neutral \emph{associates}, which extends Bonnet's classical one.}  Section~\ref{sec:conclusions} collects the conclusions of this work and a few comments on their possible bearing on the theory of plates and shells, which was our original motivation and remains our primary intent.

The paper is closed by two appendices, where we give details about proofs left out of the main text to ease the flow of our presentation. Two animations and a Maple script are also provided as supplementary material for the reader. The animations show motions traversing \nigh{bending-neutral associates of a minimal surface}: they visibly convey an impression of \emph{gliding} instead of bending, in accord with the theory. The Maple script was used to generate the animations and can be modified by the reader to explore further minimal surfaces and their deformations.

\section{Glimpses of Surface Differential Calculus}\label{sec:calculus}
To make our paper self-contained, we recall some basic notions of differential calculus on smooth surfaces in three-dimensional space, especially those playing a role in the development that follows. 

Let $\surface$ be a smooth, orientable surface in three-dimensional Euclidean space $\euclid$; it can be represented (at least locally) by a mapping of class $C^2$ defined on a domain $S\subset\mathbb{R}^2$ described by a system of coordinates. Here, however, we shall keep our discussion independent of any particular choice of coordinates.

A scalar field $\varphi:\surface\to\mathbb{R}$ is differentiable at a point $p\in\surface$, if for every curve $\curve:t\mapsto \curve(t)\in\surface$ such that $\curve(t_0)=p$ there is a vector $\nablas\varphi$ on the tangent plane  $\mathscr{T}_{p}$ to $\surface$ at $p$ such that 
\begin{equation}
	\label{eq:surface_tangent_definition}
	\left.\frac{\dd}{\dd t}\varphi(\curve(t))\right|_{t=t_0}=\nablas\varphi\cdot\tangent,
\end{equation}
where $\tangent$ is the unit tangent vector to $\curve$ at $p$. We call $\nablas\varphi$ the \emph{surface gradient} of $\varphi$.\footnote{This definition agrees with that given by Weatherburn \cite[p.\,220]{weatherburn:differential_1}, who identifies $\nablas\varphi$ with the ``vector quantity whose direction is that direction on the surface at $p$ which gives the maximum arc-rate of increase of $\varphi$, and whose magnitude is this maximum rate of increase'', which is perhaps more descriptive.}

Equivalently, $\nablas\varphi$ can be introduced by extending $\varphi$ in a three-dimensional neighborhood that contains $\surface$ to a smooth function $\tilde{\varphi}$ such that $\tilde{\varphi}|_\surface=\varphi$. Letting $\nabla\tilde{\varphi}$ denote the ordinary, three-dimensional gradient of $\tilde{\varphi}$, we can uniquely decompose it in its tangential and normal components relative to $\surface$:
\begin{equation}
	\label{eq:gradient_decomposition}
	\nabla\tilde{\varphi}=\proj\nabla\tilde{\varphi}+(\nabla\tilde{\varphi}\cdot\normal)\normal,
\end{equation}
where $\normal$ is a selected unit normal to $\surface$ at $p\in\surface$ and $\proj:=\I-\normal\otimes\normal$ is the projection onto the tangent plane $\mathscr{T}_p$. We can then set
\begin{equation}
	\label{eq:surface_normal_gradients}
	\nablas\varphi:=\proj\nabla\tilde{\varphi}\quad\text{and}\quad\nablan\tilde{\varphi}:=\nabla\tilde{\varphi}\cdot\normal,
\end{equation}
the former depending only on $\varphi$ and the latter on its extension $\tilde{\varphi}$. One can devise $\tilde{\varphi}$ so that $\nablan\tilde{\varphi}=0$, as in \cite{silhavy:new}: this is called a \emph{normal} extension of $\varphi$. Here, we shall not make use of such a restricted class of extensions. 

The advantage of defining $\nablas\varphi$ as in \eqref{eq:surface_normal_gradients} is that this definition is easily extended to higher derivatives. Denoting by the same symbol also the extension of $\varphi$ (no confusion can arise, as only the normal derivative depends on the extension), we have that 
\begin{equation}
	\label{eq:second_surface_gradient_definition}
	\nablastwo\varphi:=\nabla[\nabla\varphi-(\nabla\varphi\cdot\normal)\normal]\proj,
\end{equation}
where $\nabla$ denotes the standard gradient. By expanding \eqref{eq:second_surface_gradient_definition} and setting $\nablas\normal:=(\nabla\normal)\proj$, having also extended $\normal$ out of $\surface$, by use of \eqref{eq:surface_normal_gradients}, we arrive at 
\begin{equation}
	\label{eq:second_surface_gradient_expression}
	\nablastwo\varphi=\proj(\nabla^2\!\varphi)\proj-(\nablan\varphi)\nablas\normal-\normal\otimes\curvature(\nablas\varphi).
\end{equation}
Two consequences follow immediately from \eqref{eq:second_surface_gradient_expression}. First, $\nablastwo\varphi$ clearly depends on the extension of $\varphi$ outside $\surface$; an intrinsic definition can be achieved by resorting to normal extensions of $\varphi$. Second, since both $\nabla^2\varphi$ and the \emph{curvature tensor} $\nablas\normal$, which at each point $p\in\surface$ maps $\mathscr{T}_p$ into itself, are symmetric, the skew-symmetric part of $\nablastwo\varphi$ does not generally vanish, but it is independent of the extension of $\varphi$,
\begin{equation}
	\label{eq:skew_part_second_surface_derivative}
	\skw(\nablastwo\varphi)=\skw[\curvature(\nablas\varphi)\otimes\normal],
\end{equation}
which is both intrinsic and fully determined by $\nablas\varphi$ and the curvature of $\surface$ (see also \cite{kralj:curvature,rosso:parallel}). Letting the axial vector associated with the skew-symmetric part of $\nablastwo\varphi$ be the \emph{surface} curl of $\nablas\varphi$, we can rewrite \eqref{eq:skew_part_second_surface_derivative} as
\begin{equation}
	\label{eq:surface_curl_gradient}
	\curls\nablas\varphi=\normal\times\curvature(\nablas\varphi).
\end{equation}

More generally, for $\h$ a given \emph{tangential} vector field of class $C^1$ on $\surface$, \eqref{eq:skew_part_second_surface_derivative} is replaced by
\begin{eqnarray}
	\label{eq:integrability_vector}
	\skw(\nablas\h)=\skw(\curvature\h\otimes\normal).
\end{eqnarray} 
Similarly, for a given surface second-rank  tensor field $\F$ of class $C^1$ on $\surface$ such that $\F\normal=\zero$, 
\begin{equation}
	\label{eq:integrability_tensor}
	\skw(\nablas\F)=\skw(\F\curvature\otimes\normal),
\end{equation}
where the skew-symmetric part is taken on the last two legs of the third-rank tensors involved.

Consider now a closed curve $\curve$ over $\surface$ and denote by $\subsurface\subset\surface$ the portion of $\surface$ that it encloses. The \emph{circulation theorem} (see, for example, \cite[p.\,243]{weatherburn:differential_1}) says that for a (not necessarily tangential) vector field $\h$ defined on $\surface$ the following identity holds
\begin{equation}
	\label{eq:circulation_theorem}
	\oint_\curve\h\cdot\tangent\dd s=\int_{\subsurface}\curls\h\cdot\normal\dd A,
\end{equation}
where $\tangent$ is the unit tangent to $\curve$ (oriented anti-clockwise relative to $\normal$), $s$ is the arc-length measure along $\curve$, and $A$ is the area measure on $\surface$. 

It readily follows from \eqref{eq:circulation_theorem} and \eqref{eq:surface_curl_gradient} that a sufficient and necessary condition for a surface vector field $\h$ to be the surface gradient of a scalar field $\varphi$ is that both $\h$ and $\curls\h$ are tangential to $\surface$ \cite[p.\,244]{weatherburn:differential_2}. Said differently, for a tangential surface vector field $\h$, \eqref{eq:integrability_vector} is also a sufficient \emph{integrability} condition, which guarantees the existence (at least locally) of a scalar field $\varphi$ such that $\h=\nablas\varphi$. Similarly, \eqref{eq:integrability_tensor} is also a sufficient integrability condition for a surface second-rank tensor field $\F$ that annihilates the normal: it guarantees the existence of a surface vector field $\h$ such that $\F=\nablas\h$.

Finally, we define the \emph{surface} Laplacian as $\laps\varphi:=\divs\nablas\varphi$ and we shall say that a scalar function $\varphi$ of class $C^2$ on $\surface$ such that $\laps\varphi=0$ is \emph{surface harmonic}.\footnote{As remarked in \cite[p.\,5]{weatherburn:differential_2}, $\laps\varphi$ is just the same as the \emph{differential parameter of the second order} introduced by Beltrami (see pp.\,148, 321 of \cite{beltrami:opere}).}

\section{Bending-Neutral Deformations}\label{sec:deformations}
\begin{figure}[h]
	\centering\includegraphics[width=.66\linewidth]{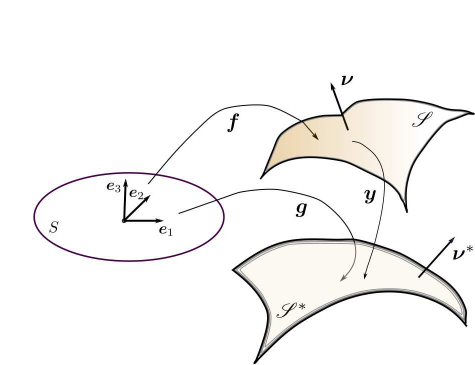}
	\caption{The surface $\surface$ is obtained from the reference configuration $S$ in the $(x_1,x_2)$ plane of a Cartesian frame $\framee$ via the deformation $\f$. $\surface^\ast$ is similarly obtained via the deformation $\g$. $\surface$ is mapped onto $\surface^\ast$ by the deformation $\y$; $\normal$ and $\normal^\ast$ are the unit vector fields that orient $\surface$ and $\surface^\ast$, respectively, while $\e_3$ is the unit normal orienting $S$.}
	\label{fig:surfaces}
\end{figure}

In this section we are concerned with the kinematics of plates. We think of $\surface$ as the image of the mid-surface of a plate $S\subset\mathbb{R}^2$ in the $(x_1,x_2)$ plane of a Cartesian frame $\framee$ under a deformation $\f:S\to\euclid$, which is assumed to be of class $C^2$ (see Fig.~\ref{fig:surfaces}).
In turn, $\surface$ can be regarded as the reference configuration of a \emph{shell}, which can be further deformed, as it will be at a later stage.

\subsection{Kinematics}\label{sec:kinematics}
Letting $\x=x_1\e_1+x_2\e_2$ denote the position vector in $S$ and extending to the present geometric setting results  well-known from three-dimensional kinematics (see, e.g., Chap.\,6 of \cite{gurtin:mechanics}), the deformation gradient $\nabla\f$ can be represented as
\begin{equation}\label{eq:representation_nabla_y}
	\nabla\f=\lambda_1\bv_1\otimes\cv_1+\lambda_2\bv_2\otimes\cv_2,
\end{equation}
where $\nabla$ denotes the gradient in $\x$, the positive scalars $\lambda_1$, $\lambda_2$ are the principal stretches and the unit vectors $\bm{c}_i$, $\bm{b}_i$ are the corresponding \emph{right} and \emph{left} principal directions of stretching. While $\cv_1(\x)$ and $\cv_2(\x)$ live on the $(\e_1,\e_2)$ plane for all $\x\in S$, $\bv_1(\x)$ and $\bv_2(\x)$ live on the tangent plane $\tplane_{\f(\x)}$ to $\surface$ at $\f(\x)$. The right and left Cauchy-Green tensors are correspondingly given by 
\begin{subequations}\label{eq:representations}
	\begin{align}
		\Cf&:=\naf\trans\naf=\lambda_1^2\cv_1\otimes\cv_1+\lambda_2^2\cv_2\otimes\cv_2,\label{eq:representation_C}\\
		\Bf&:=\naf\naf\trans=\lambda_1^2\bv_1\otimes\bv_1+\lambda_2^2\bv_2\otimes\bv_2.\label{eq:representation_B}
	\end{align}	
\end{subequations}
By the polar decomposition theorem  for the deformation of surfaces established in \cite{man:coordinate}  within the general coordinate-free theory introduced in \cite{gurtin:continuum} (see also \cite{pietraszkiewicz:determination}), the deformation gradient $\nabla\f$ can be written as
\begin{equation}
	\label{eq:polar_decomposition}
	\nabla\f=\Rf\Uf=\Vf\Rf,
\end{equation}
where $\Rf$ is a rotation of the special orthogonal group in three dimensions $\orth$ and 
\begin{equation}
	\label{eq:stretching_tensors}
	\Uf=\lambda_1\cv_1\otimes\cv_1+\lambda_2\cv_2\otimes\cv_2\quad\text{and}\quad \Vf=\lambda_1\bv_1\otimes\bv_1+\lambda_2\bv_2\otimes\bv_2
\end{equation}
are the \emph{stretching} tensors. It readily follows from \eqref{eq:polar_decomposition} that $\bv_i=\Rf\cv_i$.

\subsection{Bending and drilling rotations}\label{sec:bending}
In \eqref{eq:polar_decomposition}, the rotation $\Rf$, which can be different at different places $\x$, represents the non-metric component of the deformation $\f$. It follows from Euler's theorem on rotations (see, for example, \cite{palais:disorienting} for a modern discussion) that for every $\R\in\orth$ there is an \emph{axis}, designated by a unit vector $\e\in\sphere$, and a scalar $\alpha\in[0,\pi]$ such that $\R$ is a (conventionally anticlockwise) rotation about $\e$ by angle $\alpha$. More precisely, the following representation applies,
\begin{equation}
	\label{eq:rotation_representation}
	\R=\I+\sin\alpha\W(\e)+(1-\cos\alpha)\W(\e)^2,
\end{equation}
where $\I$ is the identity tensor, $\W(\e)$ is the skew-symmetric tensor associated with $\e$ (so that $\W(\e)\vv=\e\times\vv$ for all vectors $\vv$). $\orth$ is set by \eqref{eq:rotation_representation} into a one-to-one correspondence with the ball $\ball^3_\pi\subset\euclid$ of radius $\pi$ centred at the origin. This latter corresponds to the identity $\I$, whereas the points on $\partial\ball^3_\pi$ are $\pi$-\emph{turns}, that is, rotations by an angle $\pi$, which are the only symmetric members of $\orth$ other than $\I$.

As shown in Appendix~\ref{sec:rotation}, $\Rf$ can be uniquely decomposed into two rotations, $\Rfd$ with axis along the normal $\e_3$ to $S$ and $\Rfb$ with axis in the $(\e_1,\e_2)$ plane. We shall interpret $\Rfb$ as the \emph{bending} component of $\Rf$ and $\Rfd$ as the \emph{drilling} component.\footnote{We have retraced the first occurrence of ``drilling degrees of freedom'' to \cite{hughes:drilling}; it was soon picked up \cite{fox:drill}, see also \cite{saem:in-plane} for an  instance of more recent usage.} Formally, we write
\begin{eqnarray}
	\label{eq:rotation_decomposition}
	\Rf=\Rfb\Rfd,
\end{eqnarray}
where the order of composition matters, as bending and drilling rotations are about different axes. A pure drilling rotation, for which $\Rfb$ is uniform but $\Rfd$ is not, would pertain to a deformation $\f$ that strains $S$ but leaves it planar. A pure bending rotation, for which $\Rfd$ is uniform but $\Rfb$ is not, would instead pertain to a deformation $\f$ that brings $S$ out of a plane.

To identify both $\Rfd$ and $\Rfb$ for a given $\Rf$, it is expedient to represent a member $\R\in\orth$ via Rodrigues' formula,\footnote{The reader is referred to \cite{altmann:hamilton} for a witty account on Rodrigues' representation of rotations and its connections with Hamilton's quaternions.}
\begin{equation}
	\label{eq:rotation_vector_representation}
	\R(\av)=\frac{1}{1+a^2}\left\{(1-a^2)\I+2\av\otimes\av+2\W(\av)\right\},
\end{equation}
where $\av$ is a vector of arbitrary length $a:=\sqrt{\av\cdot\av}$. This is the \emph{vector} representation of a rotation, where $\orth$ is set into a one-to-one correspondence with the whole translation space $\translation$ of $\euclid$: the identity is still represented by the origin, but here $\pi$-turns are points at infinity.

Equation \eqref{eq:rotation_vector_representation} follows  from \eqref{eq:rotation_representation} by setting $\av=\tan(\alpha/2)\e$ and remarking that $\W(\e)^2=-\I+\e\otimes\e$. It can also be easily inverted: for given $\R\in\orth$, the representing vector $\av$ is identified by
\begin{equation}
	\label{eq:W(a)}
	\W(\av)=\frac{\R-\R\trans}{1+\tr\R},
\end{equation}
whence
\begin{equation}
	\label{eq:a_squared}
	a^2=\frac{3-\tr\R}{1+\tr\R},
\end{equation}	
which shows that a $\pi$-turn is characterized by having $\tr\R=-1$.

Rodrigues \cite{rodrigues:lois} proved a composition formula, which in modern terms can be phrased as follows (see, for example, also \cite{mladenova:vector}): for $\R(\av_1)$ and $\R(\av_2)$ in $\orth$ represented as in \eqref{eq:rotation_vector_representation},
\begin{equation}
	\label{eq:rotation_composition_formula}
	\R(\av_2)\R(\av_1)=\R(\av),
\end{equation}
where
\begin{equation}
	\label{eq:a_compostion_formula}
	\av=\frac{\av_1+\av_2+\av_2\times\av_1}{1-\av_1\cdot\av_2}.
\end{equation}
In particular, \eqref{eq:a_compostion_formula} shows that the composed rotation $\R(\av)$ is a $\pi$-turn whenever $\av_1$ and $\av_2$ belong to the cone defined by $\av_1\cdot\av_2=1$.

Let $\R(\av)$ be given in $\orth$; we prove in Appendix~\ref{sec:rotation} that it can uniquely be decomposed as in \eqref{eq:rotation_composition_formula} into a rotation $\R(\av_1)$ with $\av_1$ parallel to a given unit vector $\e$ and a rotation $\R(\av_2)$ with $\av_2$ orthogonal to $\e$. The vectors $\av_1$ and $\av_2$ are explicitly determined as
\begin{equation}
	\label{eq:a_formulas}
	\av_1=(\av\cdot\e)\e,\quad\av_2=\frac{1}{1+(\av\cdot\e)^2}\{\I+(\av\cdot\e)\W(\e)\}\proje\av,
\end{equation}
where $\proje:=\I-\e\otimes\e$ is the projector on the plane orthogonal to $\e$. It easily follows from \eqref{eq:a_formulas} that
\begin{equation}
	\label{eq:a_mod_formulas}
	a_1^2=(\av\cdot\e)^2\quad\text{and}\quad a_2^2=\frac{a^2-a_1^2}{1+a_1^2}.
\end{equation}
It is interesting to note that the rotation $\R(\av_1)$ just determined can also be characterized as the member in the subgroup $\mathsf{SO}(\e)\subset\orth$ of  all rotations with axis $\e$ that is the \emph{closest} to $\R(\av)$ in the sense made precise in Appendix~\ref{sec:rotation}, where this result is proved in detail.\footnote{Much in the spirit of the variational characterization of the polar decomposition itself \cite{martins:variational,martins:variational_1980} initiated by a kinematic note of Grioli~\cite{grioli:proprieta} recently revived and commented in \cite{neff:grioli's}.} 

For $\Rf\in\orth$, we set $\Rfd:=\R(\df)$ and $\Rfb:=\R(\bvf)$, where, repeating \eqref{eq:a_formulas} verbatim,
\begin{equation}
	\label{eq:d_and_b_reference_formulas}
	\df=(\avf\cdot\e_3)\e_3,\quad\bvf=\frac{1}{1+(\avf\cdot\e_3)^2}\{\I+(\avf\cdot\e_3)\W(\e_3)\}\mathbf{P}(\e_3)\avf,
\end{equation}
and $\avf$ is the vector associated with $\Rf$ via \eqref{eq:rotation_vector_representation}, so that \eqref{eq:rotation_decomposition} applies. \nigh{We shall call $\df$ and $\bvf$ the \emph{drilling} and \emph{bending} contents of $\Rf$, respectively. It should be noted that, as typical of any kinematic quantities, both $\df$ and $\bvf$ inherit from $\Rf$ its frame-dependent character; neither should be mistaken for \emph{measures} of drilling and bending, for which we need frame-indifferent quantities.}

\subsection{Bending neutrality}\label{eq:neutrality}
We wish to identify a class of deformations that do not alter the bending content of a surface $\surface$ (relative to its flat reference configuration $S$). Let $\surface$ be obtained as above from $S$ via the deformation $\f$ and let $\y:\surface\to\euclid$ be a deformation of $\surface$ that transforms it into $\surface^\ast$ (see Fig.~\ref{fig:surfaces}). \nigh{As for $\nabla\f$ in \eqref{eq:polar_decomposition}, we can decompose the surface gradient $\nablas\y$ of the incremental deformation $\y$ as
	\begin{equation}
		\label{eq:incremental_polar_decomposition}
		\nablas\y=\Ry\Uy,
	\end{equation}
	where $\Uy$ is the \emph{surface} stretching tensor that maps the tangent plane $\tplane_{\bm{\xi}}$ onto itself for all $\bm{\xi}\in\surface$ and $\Ry$ is a rotation of $\orth$.} We say that $\y$ is a \emph{bending-neutral} deformation if 
\begin{equation}
	\label{eq:bending_neutral_definition}
	\Ry=\nigh{\Rfd}=\Rnu\in\orthnu\quad\text{on}\ \surface,
\end{equation}
where $\normal$ is the unit normal to $\surface$ oriented coherently with the orientation of $S$. This amounts to require that the bending component $\Ryb$ of $\y$ be any \emph{uniform} rotation $\R_0$, which can always be chosen to be the identity $\I$ by an appropriate change of frame.

The surface $\surface^\ast$ can also be regarded as the image of $S$ under the deformation $\g:=\y\circ\f$. We now show that for every bending-neutral deformation $\y$ of $\surface$ the deformations $\g$ and $\f$ have one and the same bending content. To this end, we first remark that the unit normal $\normal$ at $\f(\x)$ is given by
\begin{equation}
	\label{eq:normal_field}
	\normal=\frac{\naf\e_1\times\naf\e_2}{|{\naf\e_1\times\naf\e_2}|}=\frac{\co\naf\e_3}{|\co\naf\e_3|},
\end{equation}
where $\co(\cdot)$ denotes the cofactor tensor (a definition is given, for example, in Sect.\,2.11 of \cite{gurtin:mechanics}). Since, for any two tensors $\Lt$ and $\M$, $\co(\Lt\M)=\co(\Lt)\co(\M)$ and $\co(\R)=\R$ for all $\R\in\orth$, it follows from \eqref{eq:polar_decomposition} and \eqref{eq:normal_field} that 
\begin{subequations}
	\label{eq:normal_representations}
	\begin{equation}
		\label{eq:normal_representation}
		\normal(\f(\x))=\Rf\e_3\quad\forall\ \x\in S,
	\end{equation}
	as, by \eqref{eq:stretching_tensors}, $\co(\Uf)=(\det\Uf)\e_3\otimes\e_3$. Similarly, denoting by $\normal^\ast$ the unit normal to $\surface^\ast$, we also have that 
	\begin{equation}
		\label{eq:normal_star_representation}
		\normal^\ast(\g(\x))=\Rg\e_3\quad\forall\ \x\in S.
	\end{equation}
\end{subequations}
On the other hand, since $\surface^\ast=\y(\surface)$, 
\begin{equation}
	\label{eq:normal_y_star_representation}
	\normal^\ast(\bm{\xi})=\frac{\co(\nablas\y)\normal}{|\co(\nablas\y)\normal|}\quad\forall\ \bm{\xi}\in\surface.
\end{equation}
\nigh{Moreover, since $\nablas\y$ can be decomposed as in \eqref{eq:incremental_polar_decomposition}, and $\Uy$ can be written as}
\begin{equation}
	\label{eq:surface_stretching_tensor}
	\Uy(\bm{\xi})=\mu_1\uv_1\otimes\uv_1+\mu_2\uv_2\otimes\uv_2\quad\text{with}\quad\mu_1,\mu_2>0
\end{equation}
in its eigenframe $(\uv_1,\uv_2)$ on the tangent plane $\tplane_{\bm{\xi}}$ to $\surface$, we see that $\normal^\ast=\Ry\normal$, as $\co(\Uy)=(\det\Uy)\normal\otimes\normal$. Then \eqref{eq:bending_neutral_definition} implies that 
\nigh{
	\begin{equation}\label{eq:normal_preservation}
		\normal^\ast(\g(\x))=\normal(\f(\x))\quad\text{for all}\quad\x\in S,
	\end{equation}
}
which by \eqref{eq:normal_representations} and \eqref{eq:rotation_vector_representation} amounts to the equation
\begin{equation}
	\label{eq:bending_equality}
	(\bvf-\bg)\times\e_3=\zero\quad\Rightarrow\quad \bvf=\bg=\bv\quad\text{since}\quad\bvf\cdot\e_3=\bg\cdot\e_3=0.
\end{equation}
This is the desired result, which justifies our identifying the  vector $\bv$ as the \emph{universal} bending content, as it is common to all surfaces related via a bending-neutral deformation.

\nigh{While, as already remarked, both drilling and bending contents are frame-dependent, the notion of bending-neutrality is frame-indifferent: it says that whatever bending content is assigned to a surface $\surface$ in a frame, it will be the same as the one assigned to $\surface^\ast$, if $\y$ is  bending-neutral.}

The definition of bending-neutral deformations was introduced in \cite{virga:pure} in an attempt to provide a solid foundation to the notion of \emph{pure measures} of bending in plate theory. These were defined as deformation measures (either scalars or tensors) invariant under the action of all bending-neutral deformations. For example, it was shown in \cite{virga:pure} that the tensor
\begin{equation}
	\label{eq:pure_measure}
	\pure:=(\nabla\normal)\trans(\nabla\normal)
\end{equation}
is a pure measure of bending, whereas the tensor $\mathbf{K}:=(\nabla\normal)(\nabla\normal)\trans$ is not, even if it shares with $\pure$ all scalar invariants.

\subsection{Compatibility condition}\label{sec:compatibility}
The existence of a bending-neutral deformation is subject to a compatibility condition that involves the curvature of $\surface$: as shown in \cite{virga:pure}, the rotation $\Rnu$ and the surface stretching $\Uy$ that decompose $\nablas\y$ must be such that 
\begin{equation}
	\label{eq:compatibility_condition}
	\nablas\normal=\Rnu\Uy^{-1}\curvature\Rnu\Uy.
\end{equation}
This condition ensures that the tensor $\nablas\normal^\ast$ be symmetric, as it should. We represent $\nablas\normal$ as
\begin{equation}
	\label{eq:curvature_tensor_representation}
	\nablas\normal=\kappa_1\n_1\otimes\n_1+\kappa_2\n_2\otimes\n_2,
\end{equation}
where $(\kappa_1,\kappa_2)$ are the principal curvatures of $\surface$ and $(\n_1,\n_2)$ the corresponding principal directions of curvature, and we let $\beta$ be the angle the eigenframe $(\n_1,\n_2)$ of $\nablas\normal$ makes with the eigenframe $(\uv_1,\uv_2)$ of $\Uy$ in \eqref{eq:surface_stretching_tensor}. By writing $\Rnu$ as in \eqref{eq:rotation_representation},
\begin{equation}
	\label{eq:drilling_rotation_representation}
	\Rnu(\alpha)=\I+\sin\alpha\W(\normal)-(1-\cos\alpha)\proj\quad\text{with}\quad\alpha\in[-\pi,\pi],
\end{equation}
equation \eqref{eq:compatibility_condition} can then be reduced to a single scalar equation, see \cite{virga:pure},
\begin{equation}
	\label{eq:compatibility_condition_single}
	\frac{\mu_2}{\mu_1}=1-\frac{(\kappa_1+\kappa_2)\sin\alpha}{\kappa_2\sin\alpha+(\kappa_1-\kappa_2)\cos\beta\sin(\alpha+\beta)},
\end{equation}
which is thus the general compatibility condition for the existence of a bending-neutral deformation of $\surface$. Since only the ratio between the principal surface stretches is prescribed by \eqref{eq:compatibility_condition_single}, $\Uy$ can only determined to within a multiplicative surface dilation $\lambda\proj$, with $\lambda$ an arbitrary \emph{positive} scalar field on $\surface$.

A special class of solutions of \eqref{eq:compatibility_condition_single} is obtained for $\beta=0$ and $\alpha\notin\{0,\pi\}$. This requires $\surface$ to be a \emph{hyperbolic} surface, that is, such that the Gaussian curvature $K:=\kappa_1\kappa_2<0$, as \eqref{eq:compatibility_condition_single} reduces to
\begin{equation}
	\label{eq:compatibility_condition_reduced}
	\frac{\mu_2}{\mu_1}=-\frac{\kappa_2}{\kappa_1}>0.
\end{equation}
In this class of solutions of \eqref{eq:compatibility_condition_single}, $\Uy$ can be represented as
\begin{equation}
	\label{eq:surface_stretching_single}
	\Uy=\lambda(\kappa_1\n_1\otimes\n_1-\kappa_2\n_2\otimes\n_2).
\end{equation}
It is not restrictive to assume that both $\kappa_1>0$ and $\lambda>0$.  As proved in \cite{virga:pure}, for $\Uy$ as in \eqref{eq:surface_stretching_single}, the curvature tensor of $\surface^\ast$ can be written as\footnote{It would be a simple matter for the reader to check directly that use of 
	\eqref{eq:surface_stretching_single} in  \eqref{eq:compatibility_condition} makes the latter equation identically satisfied for any $\Rnu$ as in \eqref{eq:drilling_rotation_representation} and that use of \eqref{eq:compatibility_condition_reduced} in \eqref{eq:curvature_tensor_star} makes $\nablas\normal^\ast$ correspondingly symmetric for any $\Rnu$.}
\begin{equation}
	\label{eq:curvature_tensor_star}
	\nablas\normal^\ast=\frac{\kappa_1}{\mu_1}\n_1\otimes\Rnu\n_1+\frac{\kappa_2}{\mu_2}\n_2\otimes\Rnu\n_2,
\end{equation}
which, by \eqref{eq:drilling_rotation_representation} and \eqref{eq:compatibility_condition_reduced}, implies that the total curvature of $\surface^\ast$ is
\begin{equation}
	\label{eq:total_curvature_star}
	2H^\ast:=\tr\curvaturast=\cos\alpha\left(\frac{\kappa_1}{\mu_1}+\frac{\kappa_2}{\mu_2}\right)=0.
\end{equation}
Thus, if there exists a bending-neutral deformation $\y$ of a hyperbolic surface $\surface$ with stretching tensor $\Uy$ as in \eqref{eq:surface_stretching_single}, it necessarily maps $\surface$ into a \emph{minimal} surface $\surface^\ast$.

\subsection{Integrability condition}\label{sec:integrability}
Having solved (in a special case) the compatibility condition \eqref{eq:compatibility_condition}, we found a class of possible bending-neutral deformations $\y$ of a hyperbolic surface $\surface$ characterized by
\begin{equation}
	\label{eq:bending_neutral_possible_gradient}
	\nablas\y=\Rnu(\alpha)\Uy,
\end{equation}
where $\Uy$ is given by \eqref{eq:surface_stretching_single} and both $\lambda$ and $\alpha$ are surface fields to be determined. However, there is no guarantee that \eqref{eq:bending_neutral_possible_gradient} actually represents a surface gradient. For this to be the case, the tensor field $\F:=\Rnu(\alpha)\Uy$ must satisfy the integrability condition \eqref{eq:integrability_tensor}. It is shown in Appendix~\ref{sec:integrability_appendix} that
\eqref{eq:integrability_tensor} is satisfied if and only if $\lambda$ and $\alpha$ obey the following system of equations,
\begin{equation}\label{eq:integrability_system}
	\begin{cases}
		\nablas(\lambda\kappa_1\cos\alpha)\cdot\n_2-\nablas(\lambda\kappa_2\sin\alpha)\cdot\n_1=\lambda(\kappa_1+\kappa_2)(\cos\alpha\,\con\cdot\n_1+\sin\alpha\,\con\cdot\n_2),\\
		\nablas(\lambda\kappa_1\sin\alpha)\cdot\n_2+\nablas(\lambda\kappa_2\cos\alpha)\cdot\n_1=\lambda(\kappa_1+\kappa_2)(\sin\alpha\,\con\cdot\n_1-\cos\alpha\,\con\cdot\n_2),
	\end{cases}
\end{equation}
where $\con$ is the \emph{Cartesian connector} field defined as\footnote{The definition of Cartesian connectors is further expounded in \cite{ozenda:blend,ozenda:kirchhoff}.} 
\begin{equation}
	\label{eq:connector_definition}
	\con:=(\nablas\n_1)\trans\n_2=-(\nablas\n_2)\trans\n_1.
\end{equation}
The system \eqref{eq:integrability_system} acquires a more compact form when $\surface$ itself is a \emph{minimal} surface, 
\begin{equation}
	\label{eq:integrability_minimal_surface}
	\nablas(\mu\sin\alpha)=\normal\times\nablas(\mu\cos\alpha),
\end{equation}
where $\mu:=\lambda\kappa>0$ and $\kappa:=\kappa_1=-\kappa_2$, so that \eqref{eq:bending_neutral_possible_gradient} reduces to
\begin{equation}
	\label{eq:bending_neutral_gradient}
	\nablas\y=\mu\Rnu(\alpha)\proj.
\end{equation}
It follows from \eqref{eq:integrability_minimal_surface} that 
\begin{equation}
	\label{eq:step_1}
	\begin{cases}
		\nablas(\mu\sin\alpha)\cdot\nablas(\mu\cos\alpha)=0,\\
		|\nablas(\mu\sin\alpha)|=|\nablas(\mu\cos\alpha)|.
	\end{cases}
\end{equation}
Expanding this system, we arrive at 
\begin{equation}
	\label{eq:step_3}
	\begin{cases}
		\mu\nablas\mu\cdot\nablas\alpha=0,\\
		|\nablas\mu|^2=\mu^2|\nablas\alpha|^2,
	\end{cases}
\end{equation} 
which amounts to the alternative $\nablas\mu=\pm\mu\normal\times\nablas\alpha$.
Insertion of either forms in \eqref{eq:integrability_minimal_surface} shows that the latter is only solved by $\nablas\mu=-\mu\normal\times\nablas\alpha$,
which can equivalently be written as
\begin{equation}
	\label{eq:solution}
	\nablas\alpha=\normal\times\nablas\varphi,
\end{equation}
where we have set $\varphi:=\ln\mu$.

As recalled in Sect.~\ref{sec:calculus}, for \eqref{eq:solution} to be integrable, the field $\h=\normal\times\nablas\varphi$, which is tangential to $\surface$, must also be such that $\curls\h\cdot\normal=0$. For $\varphi$ of class $C^2$, we easily see that 
\begin{equation}
	\label{eq:surface_curl_h}
	\curls\h=\curvature\nablas\varphi-(\nablastwo\varphi)\normal+(\laps\varphi)\normal-(\divs\normal)\nablas\varphi,
\end{equation}
where  $\divs\normal:=\tr\curvature=2H$ is the total curvature of $\surface$, which vanishes for a minimal surface. By use of \eqref{eq:second_surface_gradient_expression} in \eqref{eq:surface_curl_h}, we conclude that $\curls\h\cdot\normal=0$ whenever
\begin{equation}
	\label{eq:harmonic_phi}
	\laps\varphi=0.
\end{equation}
Thus, for any surface harmonic function $\varphi$ on $\surface$, we can set $\mu=\exp(\varphi)$ in \eqref{eq:bending_neutral_gradient}. Correspondingly, \eqref{eq:solution} is integrable and it too delivers a surface harmonic function $\alpha$, as it implies that, for $\alpha$ of class $C^2$,
\begin{equation}
	\label{eq:alpha_surface_laplacinan}
	\laps\alpha=\curls\normal\cdot\nablas\varphi-\normal\cdot\curls\nablas\varphi=0
\end{equation}
because $\curls\normal=\zero$, by the symmetry of the curvature tensor. This identifies (at least locally) a  drilling rotation $\Rnu$   that makes equation \eqref{eq:bending_neutral_gradient} integrable, so as to determine locally (to within a translation) a bending-neutral deformation of $\surface$. We have actually found a whole class of bending-neutral deformations of $\surface$, each represented by a surface harmonic function $\varphi$ on $\surface$.

We shall see in the following section that this bond between bending-neutral deformations of minimal surfaces and surface harmonic functions is indeed much stronger.

\section{Minimal Surfaces}\label{sec:minimal}
Here we show that every two minimal surfaces of class $C^2$ \nigh{with the \emph{same} spherical image}\footnote{\nigh{The \emph{spherical image} of a surface is the image of its normal map (or Gauss map) on $\sphere$ (see, for example, \cite[p.\,9]{dierkes:minimal}). Saying that two surfaces share the same spherical image amounts to say that they have collectively the same system of normals, which is perhaps a more colloquial expression.}} are one the image of the other under a bending-neutral deformation in the family described in Sect.~\ref{sec:deformations}. To this end, it is particularly instrumental to resort to the celebrated Enneper-Weierstrass representation of minimal surfaces (see, for example, Sect.~3.3 of \cite{dierkes:minimal}), which requires the use of \emph{isothermal} parameters $(u,v)$ to describe a surface $\surface$.

In an isothermal parameterization, the mapping $\rv:\Omega\to\euclid$ that represents $\surface$ as image of a simply connected set $\Omega$ of the $(u,v)$ plane enjoys the following properties,
\begin{equation}
	\label{eq:isothermal_parameterization}
	|\rv_{,u}|=|\rv_{,v}|\quad\text{and}\quad\rv_{,u}\cdot\rv_{,v}=0,
\end{equation}
where a comma denotes differentiation.\footnote{The existence of isothermal parameters for a generic surface is not obvious, and for surfaces of class $C^1$ it is even not generally true. For surfaces of class $C^2$, they always exist (at least locally) by a classical theorem (see also \cite{chern:elementary}), whose proof is indeed much easier for minimal surfaces (see, for example, \cite[p.\,31]{osserman:survey}).} Moreover, since also 
\begin{equation}
	\label{eq:laplacian_isothermal}
	\lap\rv=\rv_{,uu}+\rv_{,vv}=\zero
\end{equation}
(see, for example, \cite[p.\,268]{spivack:comprehensive_4}), for any scalar-valued function $\varphi$ defined on $\surface$ in terms of the $(u,v)$ parameters,
\begin{equation}
	\label{eq:surface_laplacian_isothermal}
	\laps\varphi=\frac{1}{|\rv_{,u}|^2}\lap\varphi=\frac{1}{|\rv_{,u}|^2}(\varphi_{,uu}+\varphi_{,vv}),
\end{equation}
so that $\varphi$ is surface harmonic on $\surface$ whenever it is harmonic in $\Omega$.

\subsection{Weierstrass representation}\label{sec:Weierstrass}
Let $w:=u+iv\in\mathbb{C}$ be a complex representation of the $(u,v)$ plane. A minimal surface $\surface$ such that the normal $\normal$ regarded as a map from $\Omega$ into the unit sphere $\sphere$ is injective, but does not cover the whole of $\sphere$, can be represented (to within a translation) by a \emph{holomorphic} function $F:\Omega\to\mathbb{C}$ as
\begin{equation}
	\label{eq:Weierstrass_representation}
	\rv(w)=\Re\left(\frac12\int(1-w^2)F(w)\dd w\,\e_1+\frac{i}{2}\int(1+w^2)F(w)\dd w\,\e_2+\int wF(w)\dd w\,\e_3\right),
\end{equation}
where $\framee$ is a Cartesian frame in $\translation$ and $\Re$ denotes the real part (see, for example, \cite[p.\,117]{dierkes:minimal}). $F$ is also called the \emph{Weierstrass function} that represents $\surface$ (see also \cite{weierstrass:untersuchungen_berlin,weierstrass:untersuchungen_cambridge}).

A number of interesting consequences follow from \eqref{eq:Weierstrass_representation}, which can easily be interpreted geometrically.

We first write $F$ as
\begin{equation}
	\label{eq:F_representation}
	F=\ex^{\Phi+i\chi}=\ex^\Phi(\cos\chi+i\sin\chi).
\end{equation}
Since $F$ is holomorphic, the Cauchy-Riemann relations imply that both $\Phi$ and $\chi$ are harmonic in $\Omega$. Explicit computations give
\begin{subequations}\label{eq:Weierstrass_metric_vectors}
	\begin{align}
		\rv_{,u}&=\ex^\Phi\Big\{\Big[uv\sin\chi+\frac12(1-u^2+v^2)\cos\chi\Big]\e_1-\Big[uv\cos\chi+\frac12(1+u^2-v^2)\sin\chi\Big]\e_2\nonumber\\&+(u\cos\chi-v\sin\chi)\e_3\Big\},\label{eq:Weierstrass_metric_vector_u}\\
		\rv_{,v}&=\ex^\Phi\Big\{\Big[uv\cos\chi-\frac12(1-u^2+v^2)\sin\chi\Big]\e_1+\Big[uv\sin\chi-\frac12(1+u^2-v^2)\cos\chi\Big]\e_2\nonumber\\&-(u\sin\chi+v\cos\chi)\e_3\Big\},\label{eq:Weierstrass_metric_vector_v}
	\end{align}
	which satisfy
	\begin{equation}\label{eq:Weierstrass_metric_vector_u_v}
		|\rv_{,u}|=|\rv_{,v}|=\frac12\ex^\Phi(u^2+v^2+1)\quad\text{and}\quad\rv_{,u}\cdot\rv_{,v}=0.
	\end{equation}
\end{subequations}
Letting $\e_u$ and $\e_v$ denote the unit vectors of $\rv_{,u}$ and $\rv_{,v}$, respectively, we form an orthonormal frame $(\e_u,\e_v)$ on the tangent plane $\tplane_{\rv(w)}$ on $\surface$; the unit normal $\normal$ is then obtained from
\begin{equation}
	\label{eq:Weierstrass_normal}
	\normal=\e_u\times\e_v=\frac{1}{u^2+v^2+1}\{2u\e_1+2v\e_2+(u^2+v^2-1)\e_3\}.
\end{equation}
\nigh{This equation shows clearly the geometrical meaning of the $(u,v)$ coordinates in Weierstrass' representation: $w$ is the stereographic projection of $\normal$ (on the equatorial plane of $\sphere$ with the North pole identified with $\e_3$). Thus, all functions $F$ defined on the \emph{same} domain $\Omega$ represent minimal surfaces with the \emph{same} spherical image.}

Consider now a trajectory in $\Omega$ parameterized by $(u(t),v(t))$; $\rv$  maps it into a trajectory on $\surface$ for which
\begin{equation}
	\label{eq:Weierstrass_r_dot}
	\dot{\rv}=|\rv_{,u}|(\dot{u}\e_u+\dot{v}\e_v),
\end{equation}
where a superimposed dot denotes differentiation with respect  to $t$. Correspondingly,
\begin{equation}
	\label{eq:Weierstrass_normal_dot}
	\dot{\normal}=\dot{u}\normal_{,u}+\dot{v}\normal_{,v}=\curvature\dot{\rv}.
\end{equation}
By differentiating with respect to $u$ and $v$ both sides of \eqref{eq:Weierstrass_normal} and projecting both $\normal_{,u}$ and $\normal_{,v}$ in the frame $(\e_u,\e_v)$, we obtain that 
\begin{subequations}\label{eq:Weierstrass_partial_normal}
	\begin{align}
		\normal_{,u}&=\frac{2}{u^2+v^2+1}(\cos\chi\e_u-\sin\chi\e_v),	\label{eq:Weierstrass_partial_normal_u}\\
		\normal_{,v}&=-\frac{2}{u^2+v^2+1}(\sin\chi\e_u+\cos\chi\e_v).\label{eq:Weierstrass_partial_normal_v}
	\end{align}
\end{subequations}
Inserting \eqref{eq:Weierstrass_partial_normal} and \eqref{eq:Weierstrass_normal} into \eqref{eq:Weierstrass_normal_dot}, with the aid of \eqref{eq:Weierstrass_metric_vector_u_v}, since $\dot{u}$ and $\dot{v}$ are arbitrary, we arrive at the following representation of the curvature tensor in terms of the functions $\Phi$ and $\chi$,
\begin{equation}
	\label{eq:Weierstrass_curvature_tensor}
	\nablas\normal=\frac{4}{\ex^\Phi(u^2+v^2+1)^2}\Rnu(-\chi)(\e_u\otimes\e_u-\e_v\otimes\e_v),
\end{equation}
from which we obtain at once that
\begin{equation}
	\label{eq:Weierstrass_curvatures}
	2H=\tr\curvature=0\quad\text{and}\quad K=\det\curvature=-\frac{16}{\ex^{2\Phi}(u^2+v^2+1)^4}.
\end{equation}
Comparing \eqref{eq:Weierstrass_curvature_tensor} and \eqref{eq:curvature_tensor_representation}, we easily see that 
\begin{equation}
	\label{eq:Weierstrass_principal_curvatures}
	\kappa_1=-\kappa_2=\frac{4}{\ex^\Phi(u^2+v^2+1)^2},\quad\n_1=\cos\frac{\chi}{2}\e_u-\sin\frac{\chi}{2}\e_v,\quad \n_2=\sin\frac{\chi}{2}\e_u+\cos\frac{\chi}{2}\e_v.
\end{equation}

\subsection{Bending-neutral associates}\label{sec:associates}
Let $\surface^\ast$ be a minimal surface, different from $\surface$ \nigh{but with the same spherical image}, represented by the holomorphic function 
\begin{equation}
	\label{eq:F_star_representation}
	F^\ast=\ex^{\Phi^\ast+i\chi^\ast}=\ex^{\Phi^\ast}(\cos\chi^\ast+i\sin\chi^\ast),
\end{equation}
where both $\Phi^\ast$ and $\chi^\ast$ are harmonic functions in $\Omega$. Precisely as in \eqref{eq:Weierstrass_representation}, $F^\ast$ induces the mapping $\rv^\ast:\Omega\to\euclid$ that represents $\surface^\ast$. The vectors $\rv^\ast_{,u}$ and $\rv^\ast_{,v}$ are still given by \eqref{eq:Weierstrass_metric_vectors} with $\Phi$ replaced by $\Phi^\ast$ and $\chi$ by $\chi^\ast$. The unit vectors $(\e_u^\ast,\e_v^\ast)$ that identify a frame on the tangent plane $\tplane_{\rv^\ast(w)}$ to $\surface^\ast$ are defined accordingly and the unit normal $\normal^\ast=\e_u^\ast\times\e_v^\ast$ is still given by \eqref{eq:Weierstrass_normal} in the Cartesian frame $\framee$, so that $\normal^\ast\equiv\normal$, although $\rv^\ast\neq\rv$. It is not difficult to show by direct computation that 
\begin{equation}
	\label{eq:Weierstrass_starred_frame}
	\e_u^\ast\cdot\e_u=\e_v^\ast\cdot\e_v=\cos(\chi^\ast-\chi).
\end{equation}

Let $\y$ be a deformation that maps $\surface$ onto $\surface^\ast$. For any trajectory $(u(t),v(t))$ in $\Omega$, we have that 
\begin{equation}
	\label{eq:Weierstrass_r_star_dot}
	\dot{\rv}^\ast=\nay\dot{\rv},
\end{equation}
which by \eqref{eq:Weierstrass_metric_vector_u_v} and \eqref{eq:Weierstrass_r_dot}, along with the corresponding equations for $\rv^\ast$, can also be written as
\begin{equation}
	\label{eq:Weierstrass_r_star_dot_rewtitten}
	\ex^{\Phi^\ast}(\dot{u}\e_u^\ast+\dot{v}\e_v^\ast)=\ex^{\Phi}(\dot{u}\e_u+\dot{v}\e_v).
\end{equation}
Since the frames $(\e_u^\ast,\e_v^\ast)$ and $(\e_u,\e_v)$ are on parallel planes (both orthogonal to $\normal$) and both $\dot{u}$ and $\dot{v}$ are arbitrary, \eqref{eq:Weierstrass_r_star_dot_rewtitten} and \eqref{eq:Weierstrass_starred_frame} imply  that
\begin{equation}
	\label{eq:Weierstrass_nabla_y}	
	\nablas\y=\ex^{\Phi^\ast-\Phi}\Rnu(\chi^\ast-\chi)\proj.
\end{equation}
By comparing \eqref{eq:Weierstrass_nabla_y} and \eqref{eq:bending_neutral_gradient}, we readily conclude that $\y$ is a bending-neutral deformation with
\begin{equation}
	\label{eq:mu_alpha_identification}
	\mu=\ex^\varphi=\ex^{\Phi^\ast-\Phi}\quad\text{and}\quad\alpha=\chi^\ast-\chi
\end{equation}
directly derived from the Weierstrass functions representing the surfaces $\surface^\ast$ and $\surface$. This proves our claim that \emph{any two minimal surfaces \nigh{with the same spherical image} are  associated through a bending-neutral deformation}.

In the classical theory of minimal surfaces a special role is played by the families of \emph{associate surfaces}. Given a minimal surface $\surface$, the classical family of its associates $\surface_\theta=\y_\theta(\surface)$, described by a real parameter $\theta$, consists of minimal surfaces in isometric correspondence with one another  and such that \nigh{normals are preserved} for all $\theta$. The mapping  that changes $\surface$ into $\surface_\theta$ is also called \emph{Bonnet's transformation} (see, for example, Sect.\,3.1 of \cite{dierkes:minimal}); in the Weierstrass representation, it has the following simple form \cite[p.\,119]{dierkes:minimal},
\begin{equation}
	\label{eq:Weierstrass_Bonnet_transformation}
	F(w)\mapsto F_\theta(w):=\ex^{i\theta}F(w)
\end{equation}
with $\theta$ a constant.

Perhaps the best-known of these families is that transforming catenoids  into helicoids, for which
\begin{equation}
	\label{eq:Weierstrass_cathenoid_helicoid}
	F_\theta=\frac{c}{w^2}(\cos\theta+i\sin\theta),
\end{equation}
where $c\in\mathbb{R}$ is a constant and $\theta\in[0,\frac{\pi}{2}]$. For $\theta=0$ or $\theta=\frac{\pi}{2}$, $F_\theta$ represents a catenoid or a helicoid, respectively; for $0<\theta<\frac{\pi}{2}$, the intermediate associate surface $\surface_\theta$ is \emph{Scherk's second surface} (see, for example, \cite[p.\,148]{dierkes:minimal}).

Bending-neutral deformations  generate  far more general families of associate minimal surfaces through the transformation
\begin{equation}
	\label{eq:Weierstrass_generalized_Bonnet_transformation}
	F(w)\mapsto F^\ast(w):=\ex^{\varphi+i\vt}F(w),
\end{equation}
where both $\vp$ and $\vt$ are harmonic functions. Furthermore, in this broad sense, any minimal surface is the \nigh{bending-neutral} associate of another \nigh{with the same spherical image}. Like Bonnet's, the transformation \eqref{eq:Weierstrass_generalized_Bonnet_transformation} preserves \nigh{normals}, and it is isometric for $\vp\equiv0$. In particular, by \eqref{eq:mu_alpha_identification}, we may say that for any $\theta\in(0,\frac{\pi}{2})$ Scherk's second surface $\surface_\theta$ is the image of a catenoid (or a helicoid) under an isometric bending-neutral deformation with $\mu\equiv1$ and $\alpha\equiv\theta$.

\subsection{Universal bending content}\label{sec:universal}
We learned in Sect.~\ref{sec:deformations} that in general two surfaces related by a bending-neutral deformation have one and the same bending content. \nigh{Here, we prove that the bending content $\bv$ of a minimal surface is as universal as the normal $\normal$ in Weierstrass' representation.}

Interpreting $\rv$ as a deformation from a reference configuration $\Omega$, by \eqref{eq:polar_decomposition}, we can write $\nabla\rv=\R_{\rv}\U_{\rv}$, where
\begin{equation}
	\label{eq:U_r_and_R_r}
	\U_{\rv}=\frac12\ex^\Phi(u^2+v^2+1)\Proj\quad\text{and}\quad\R_{\rv}=\e_u\otimes\e_1+\e_v\otimes\e_2+\normal\otimes\e_3,
\end{equation}
having chosen $(\e_1,\e_2)$ as the frame of the $(u,v)$ coordinates in $\Omega$ with $\e_3=\e_1\times\e_2$. By use of \eqref{eq:W(a)}, we identify the vector $\av_{\rv}$ that represents $\R_{\rv}$ in \eqref{eq:rotation_vector_representation},
\begin{subequations}
	\begin{equation}
		\label{eq:a_r_representation}
		\av_{\rv}=\frac{1}{u\sin\chi+v(1+\cos\chi)}\{-(1+\cos\chi)\e_1+\sin\chi\e_2+[v\sin\chi-u(1+\cos\chi)]\e_3\}.
	\end{equation}
	Then we obtain from \eqref{eq:d_and_b_reference_formulas} the vectors $\dv_{\rv}$ and $\bv_{\rv}$ that represent the drilling and bending components $\R_{\rv}^{(\mathrm{b})}$ and $\R_{\rv}^{(\mathrm{d})}$ of $\R_{\rv}$,
	\begin{align}
		\dv_{\rv}&=\frac{v\sin\chi-u(1+\cos\chi)}{u\sin\chi+v(1+\cos\chi)}\e_3=-\cot\left(\frac{\chi}{2}+\phi\right)\e_3,\label{eq:d_r}\\
		\bv_{\rv}&=\frac{1}{u^2+v^2}(-v\e_1+u\e_2)=\frac{1}{\rho}\e_\phi,\label{eq:b_r}
	\end{align}
\end{subequations}
where we have expressed $(u,v)$ in polar form,
\begin{equation}
	\label{eq:polar_form}
	u=\rho\cos\phi,\quad v=\rho\sin\phi,
\end{equation}
and set $\e_\phi:=-\sin\phi\,\e_1+\cos\phi\,\e_2$.
It is clear from \eqref{eq:b_r} that $\bv_{\rv}=\bv$ is the same for all minimal surfaces, as it is independent of both $\Phi$ and $\chi$, \nigh{precisely like $\normal$ in \eqref{eq:Weierstrass_normal}.}

In a similar way, we can prove directly that the pure measure of bending $\pure$ in \eqref{eq:pure_measure} is the same for all minimal surfaces. Since 
\begin{equation}\label{Weierstrass_nabla_normal}
	\nabla\normal=\normal_{,u}\otimes\e_1+\normal_{,v}\otimes\e_2,
\end{equation}
by use of \eqref{eq:Weierstrass_partial_normal}, we give $\pure$ the following explicit form,
\begin{equation}
	\label{eq:Weierstrass_pure}
	\pure=\frac{4}{(u^2+v^2+1)^2}\Proj.
\end{equation}

\subsection{Illustration by Bour's surfaces}\label{sec:Bour}
To illustrate graphically our findings, we consider an exemplary family $\surface_t$ of bending-neutral associates of two minimal surfaces $\surface$ and $\surface^\ast$, both with historical pedigree. The parameter $t$  ranges in the interval $[0,1]$, so that $\surface_t$ is $\surface$ for $t=0$, and $\surface^\ast$ for $t=1$. The Weierstrass function representing $\surface_t$ is
\begin{equation}
	\label{eq:F_t}
	F_t(w)=w^t=\rho^t\ex^{it\phi},
\end{equation}
where we have set $w=\rho\exp(i\phi)$. It is immediately seen that $\surface_t$ is the image of $\surface$ under the bending-neutral deformation characterized by $\mu=\rho^t$ and $\alpha=t\phi$.

Here $\surface$ is Enneper's surface, while $\surface^\ast$ is Bour's surface of index $m=3$. The Weierstrass function of the former is $F(w)\equiv1$, while that of the latter is obtained by setting $m=3$ into the general formula
\begin{equation}
	\label{eq:Weierstrass_Bour_surfaces}
	F_\mathrm{B}(w):=cw^{m-2},
\end{equation}
where $c\neq0$ is a complex constant and $m\in\mathbb{R}$ (see \cite[p.\,156]{dierkes:minimal}).\footnote{Bour~\cite{bour:theorie} proved that the surfaces represented by \eqref{eq:Weierstrass_Bour_surfaces} constitute the complete set of minimal surfaces that can be mapped isometrically onto a surface of revolution. The paper \cite{bour:theorie} is the published version of the manuscript that won Bour the 1859 \emph{Grand Prix des Math\'ematiques} granted by the French Academy of Sciences in Paris. A terse account on Bour's manuscript and two others competing with it for the same prize can be found in \cite{cogliati:origins}.} Actually, $\surface_t$ is a Bour's surface for all $t\in[0,1]$, with $m=t+2$.

\begin{figure}[]
	\centering
	\begin{subfigure}[c]{0.36\textwidth}
		\centering
		\includegraphics[width=\textwidth]{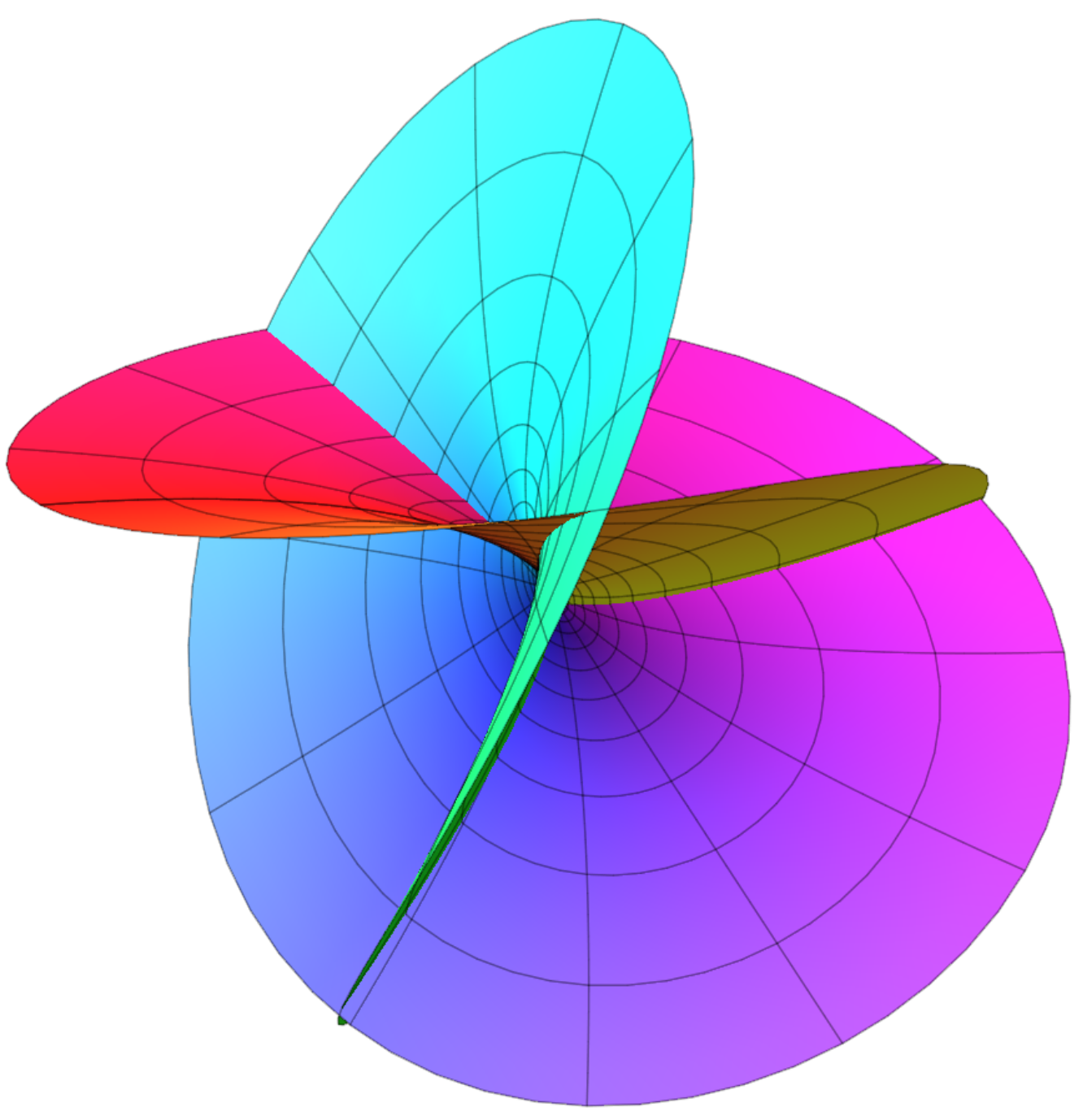}
		\caption{$t=0$}
		\label{fig:a}
	\end{subfigure}
	\hfill\\
	\begin{subfigure}[c]{0.30\textwidth}
		\centering
		\includegraphics[width=\textwidth]{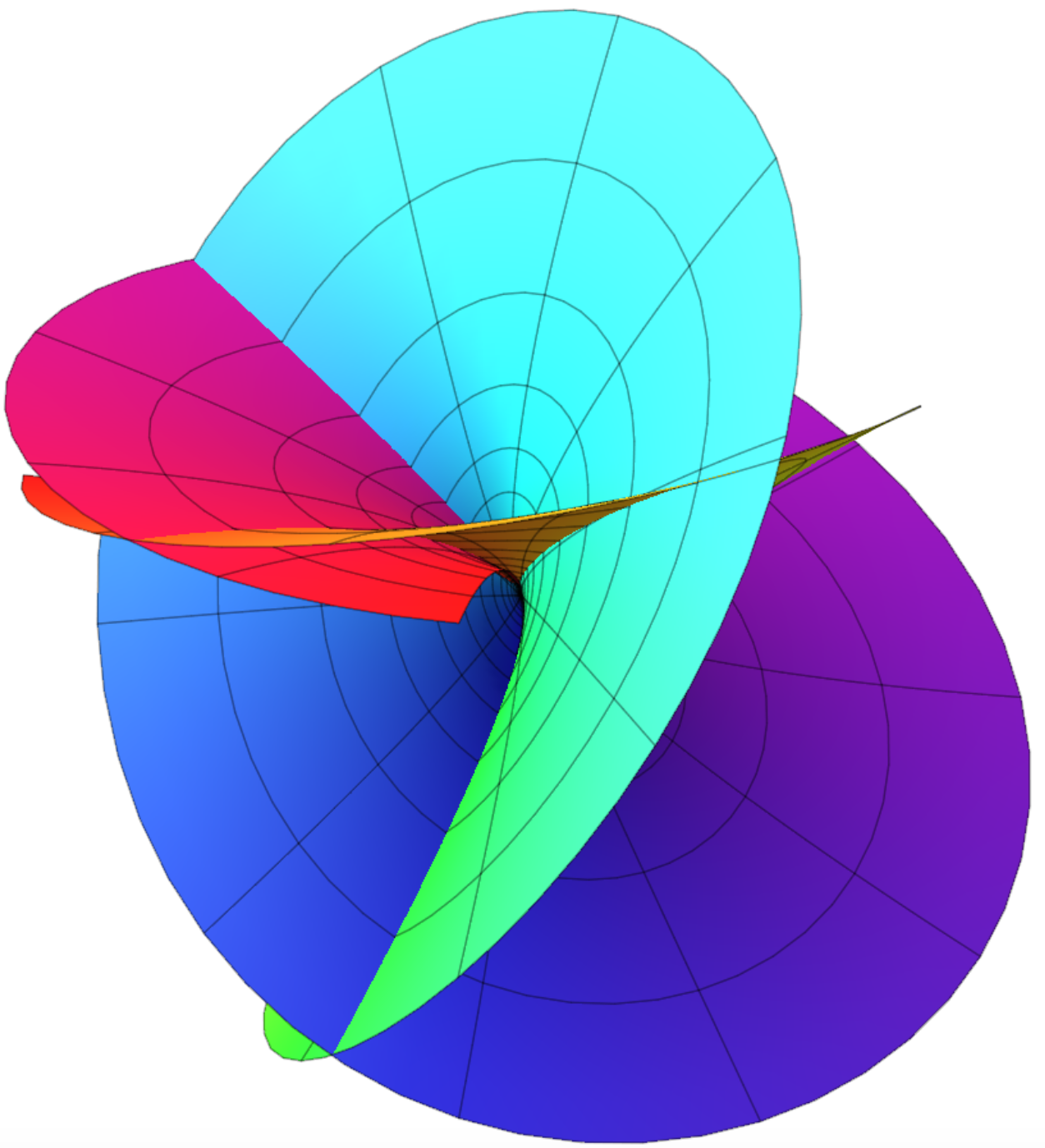}
		\caption{$t=\frac14$}
		\label{fig:b}
	\end{subfigure}
	\hfill
	\begin{subfigure}[c]{0.32\textwidth}
		\centering
		\includegraphics[width=\textwidth]{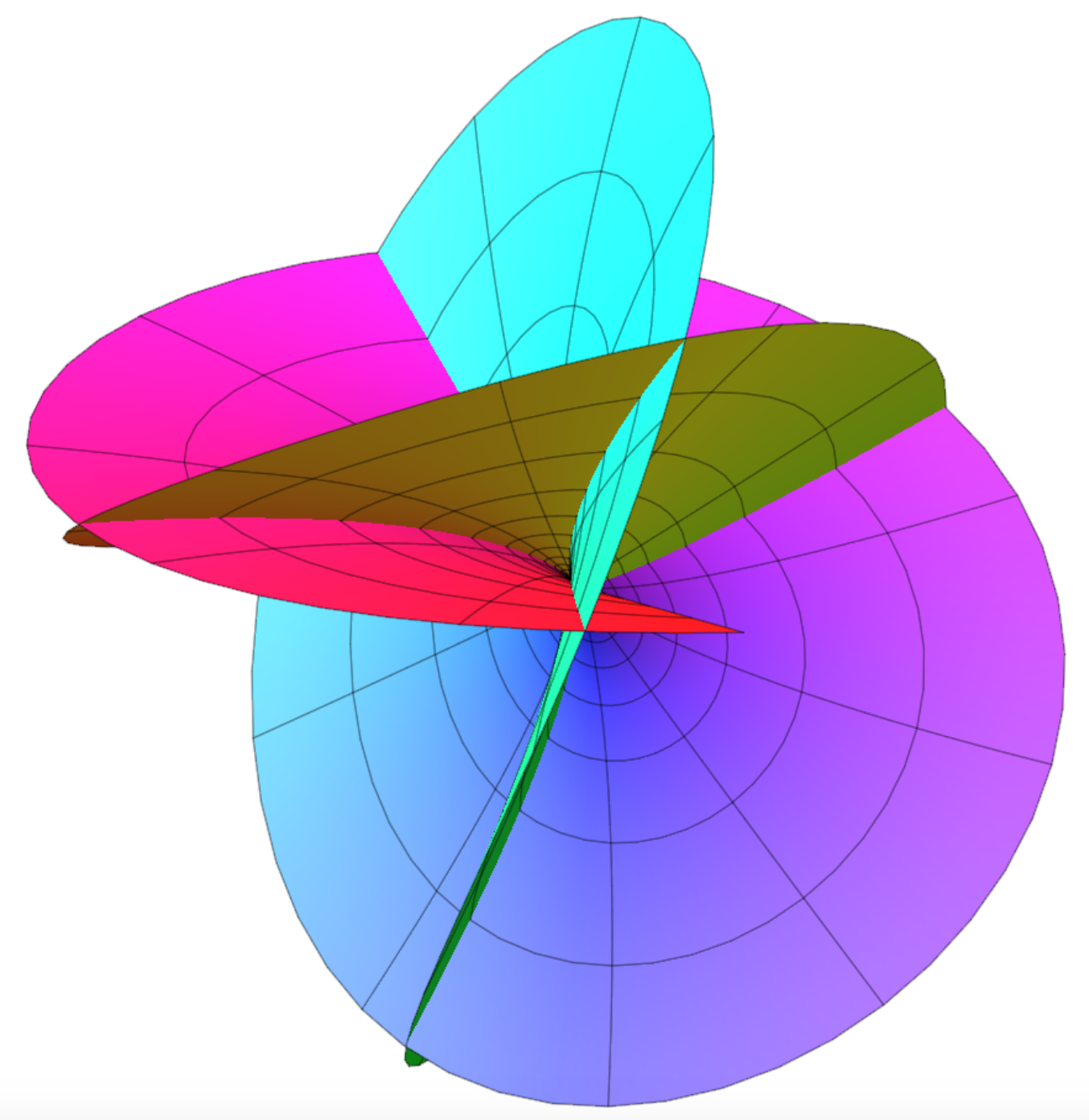}
		\caption{$t=\frac12$}
		\label{fig:c}
	\end{subfigure}
	\hfill
	\begin{subfigure}[c]{0.32\textwidth}
		\centering
		\includegraphics[width=\textwidth]{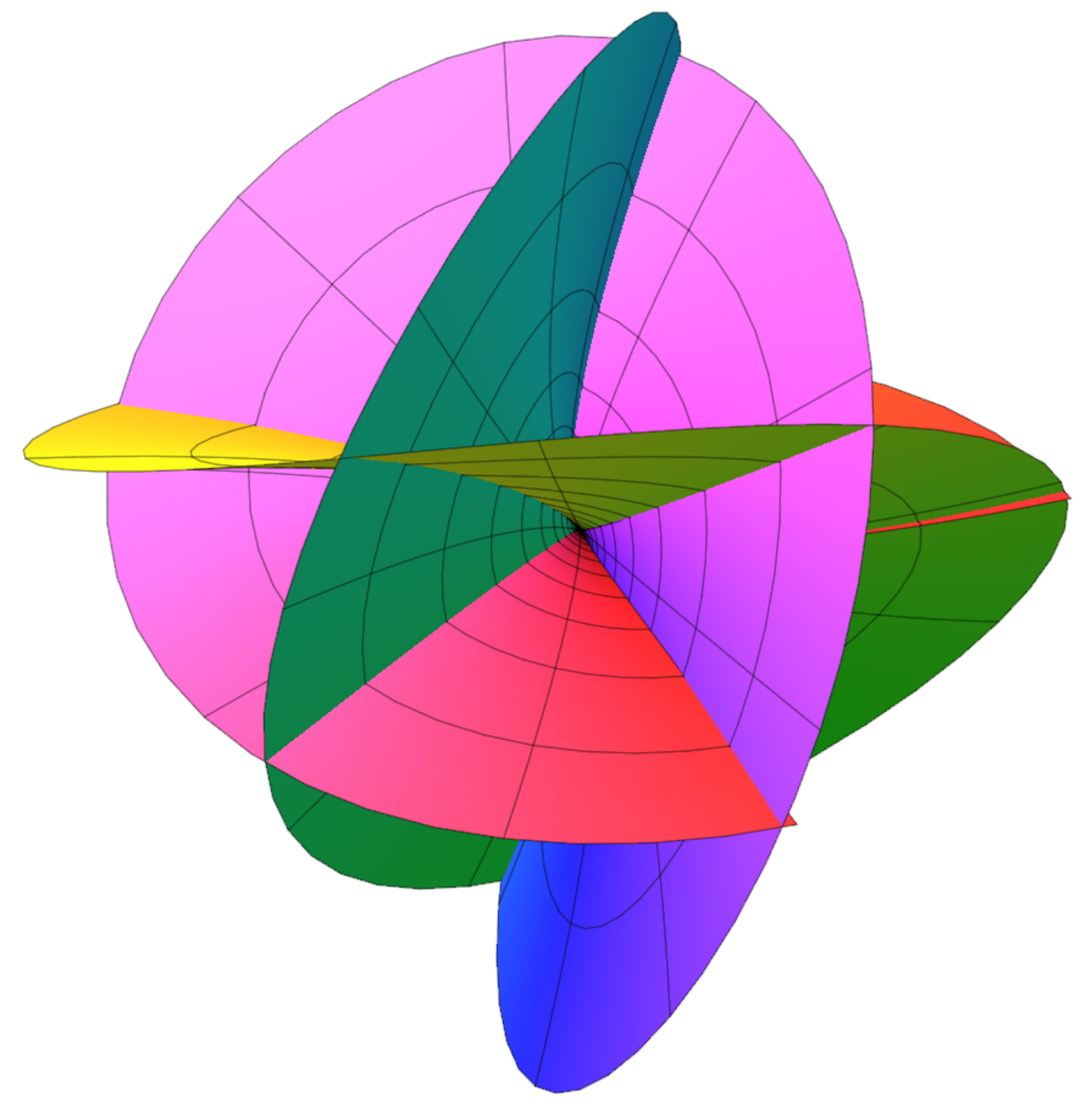}
		\caption{$t=\frac34$}
		\label{fig:d}
	\end{subfigure}
	\\
	\begin{subfigure}[c]{0.36\textwidth}
		\centering
		\includegraphics[width=\textwidth]{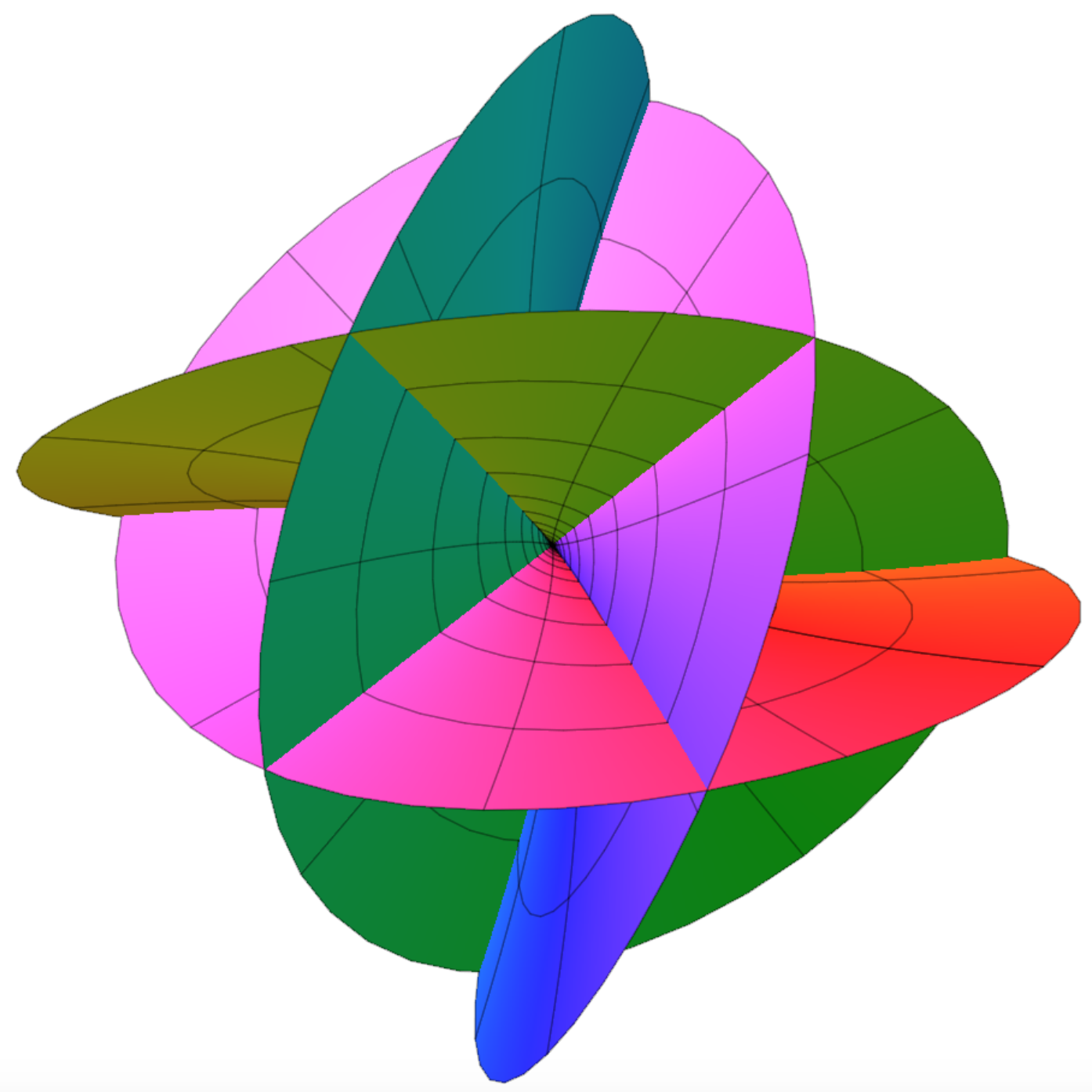}
		\caption{$t=1$}
		\label{fig:e}
	\end{subfigure}
	\hfill
	\caption{Surfaces $\surface_t$ extracted from the family described by the Weierstrass function $F_t$ in \eqref{eq:F_t}. They are images via the mapping $\rv_t(w)$ induced by $F_t(w)$ through \eqref{eq:Weierstrass_representation} of the disk $\ball_R$ of radius $R$ centred at $w=0$; the scaling factor of each image is $|\rv_t(R)|$.}
	\label{fig:gallery}
\end{figure}

Figure~\ref{fig:gallery} shows a gallery of surfaces extracted from the family $\surface_t$ represented by \eqref{eq:Weierstrass_Bour_surfaces}: these are images of a disk $\ball_R$ centred at $w=0$ (described in polar coordinates by $0<\rho<R$ and $-\pi\leq\phi\leq\pi$); the colour coding is such that each individual radius of $\ball_R$ conveys on $\surface_t$ a specific colour.

As a supplement to this paper, we also provide an animation showing the whole family $\surface_t$. We think that this makes it apparent how the moving surface progressively unfolds \emph{gliding} over itself, with \emph{no} sign of \emph{bending}. In a companion animation, we also capture the motion that illustrates the classical association between catenoid and helicoid described by \eqref{eq:Weierstrass_cathenoid_helicoid}: it clearly conveys the same gliding impression.

We finally remark that by the change of variable $w\mapsto1/w$ in \eqref{eq:Weierstrass_representation} for $F=F_\mathrm{B}$ as in \eqref{eq:Weierstrass_Bour_surfaces} Bour's surfaces are shown to be invariant (to within a reflection) under the transformation $m\mapsto-m$. This makes catenoids and helicoids somehow special members of Bour's family, as they are obtained from \eqref{eq:Weierstrass_Bour_surfaces} for $m=0$.

\section{Conclusions}\label{sec:conclusions}
Understanding how surfaces deform is the proper foundation for a sound, intrinsic  theory of plates and shells (as abundantly discussed, for example, in Sect.~XIV.13 of \cite{antman:nonlinear}). Expressing the appropriate energy contents involved in the elasticity of these bodies requires a clear identification of the local independent modes for changing their shape.

Intuitively, we may say that there are only two such modes, that is, \emph{stretching} and \emph{bending}. Rigorously, we may associate unambiguously the former mode with the stretching component $\U$ of the surface deformation gradient $\F$, as delivered by the polar decomposition theorem. But associating the bending mode with the whole orthogonal component $\R$ of $\F$ cannot be right, as this deformation \nigh{descriptor} would in general fail to be intrinsic to the surface.

By decomposing uniquely $\R$ into a \emph{drilling} rotation $\R^{\mathrm{(d)}}$ about the surface normal $\normal$ and a complementary rotation $\R^{(\mathrm{b})}$ about an axis orthogonal to $\normal$, we identified the latter as the correct \emph{bending} descriptor.
\emph{Bending-neutral} deformations were then naturally defined as those 
deformations that would leave the bending content unchanged, while affecting the drilling one. \nigh{These deformations preserve normals.} In a variational elastic theory for shells, these deformations, if existent, would not affect the bending energy, while still prompting effective changes in shape, at the expenses of a drilling energy.

Whether bending-neutral deformations do actually exist is not a trivial issue. In this paper, we gave a necessary and sufficient integrability condition that guarantees their existence, al least locally. A large class of surfaces can indeed be subjected to bending-neutral deformations:  \emph{all} minimal surfaces. Actually, there is more to it: (1) every minimal surface is transformed into a minimal surface by a bending-neutral deformation, (2) given two minimal surfaces \nigh{with the same system of normals}, there is a bending-neutral deformation that maps one into the other, \nigh{and (3) \emph{all}  minimal surfaces have indeed a \emph{universal} bending content}.

All  conceivable \nigh{normal preserving} motions driving one minimal surface into another \nigh{can only consist of} drilling and stretching, as is apparent from the gliding motions shown in the animations accompanying this paper.

In the light of our finer decomposition of surface deformations, speaking of ``a bending procedure which, at every stage, passes through a minimal surface'' \cite[p.\,102]{dierkes:minimal}, as done in almost every textbook, would be inaccurate \nigh{whenever normals are preserved}: there is indeed \emph{no} way to \emph{bend} a minimal surface, \nigh{other than by altering the system of its normals.} 

Although we deem our kinematic analysis of surfaces possibly also relevant to the energetics of elastic plates and shells, we have not attempted  any systematic treatment of this issue. \nigh{A preliminary step toward this aim would be to extract from both $\df$ and $\bvf$ objective measures of drilling and bending that would feature in a stored-energy functional.}

We are also guilty of another omission. Our deformation study was local in nature: we  left out all global considerations. This came with the price of allowing self-intersecting surfaces as deformed images, which is forbidden  in continuum mechanics by the principle of identification of material body points.

\appendix

\section{Rotation Decomposition}\label{sec:rotation}
Here we use equation \eqref{eq:rotation_composition_formula} for the decomposition of a rotation $\R(\av)\in\orth$ to prove that the vectors $\av_1$ and $\av_2$ can be identified uniquely by requiring the former to be parallel to a given unit vector $\e$ and the latter orthogonal to it.

Let $\av\in\translation$ and $\e\in\sphere$ be given so that $\av\times\e\neq\zero$.\footnote{Were $\av\times\e=\zero$, our claim would be trivially valid with $\av_1=\av$ and $\av_2=\zero$.} We solve \eqref{eq:a_compostion_formula} with $\av_1=u\e$ and $\av_2\cdot\e=0$. Projecting both sides of \eqref{eq:a_compostion_formula} along $\e$ and on the plane orthogonal to it, we find that 

\begin{subequations}\label{eq:projections}
	\begin{align}
		u&=\av\cdot\e,\label{eq:projection_e}\\
		\proje\av&=\av_2-u\e\times\av_2.\label{eq:projection_Pe}
	\end{align}
\end{subequations}
Representing $\av_2$  as
\begin{equation}
	\label{eq:a_2_representation}
	\av_2=b_1\proje\av+b_2\e\times\proje\av,
\end{equation}
where $b_1$ and $b_2$ are scalar parameters to be determined, we easily see that the only solution of \eqref{eq:projection_Pe} is obtained for
\begin{equation}
	\label{eq:b_solutions}
	b_1=\frac{1}{1+u^2}\quad\text{and}\quad b_2=\frac{u}{1+u^2}.
\end{equation}
Use of \eqref{eq:b_solutions} and \eqref{eq:projection_e} in \eqref{eq:a_2_representation} readily leads us to \eqref{eq:a_formulas} in the main text.

To prove that $\R(\av_1)$ with $\av_1$ as in \eqref{eq:a_formulas} is indeed the rotation of $\mathsf{SO}(\e)$ closest to $\R(\av)$, we introduce the (squared) distance
\begin{equation}
	\label{eq:rotation_distance}
	d(u):=|\R(\av)-\R(u\e)|^2:=\tr\{[\R(\av)-\R(u\e)]\trans[\R(\av)-\R(u\e)]\}.
\end{equation}
A direct computation shows that 
\begin{equation}
	\label{eq:d(u)}
	d(u)=6+\frac{2}{1+a^2}\frac{1}{1+u^2}\{[1+a^2-4(\av\cdot\e)^2]u^2-8(\av\cdot\e)u+a^2-3\}
\end{equation}
and it is then a simple matter to conclude that $d$ attains its minimum for $u=\av\cdot\e$ and its maximum for $u=-1/\av\cdot\e$.

\section{Gradient Integrability}\label{sec:integrability_appendix}
In this appendix we gather a number of details needed to derive the integrability conditions in \eqref{eq:integrability_system}. The objective is to enforce \eqref{eq:integrability_tensor} for the surface tensor field $\F=\Rnu(\alpha)\Uy$ where $\Uy$ is given by \eqref{eq:surface_stretching_single} and $\Rnu(\alpha)$ by \eqref{eq:drilling_rotation_representation}.

We begin by writing $\F$ explicitly in the eigenframe $\framen$ of the curvature tensor,
\begin{equation}
	\label{eq:F_curvature_frame_representation}
	\F=\lambda\{\kappa_1\cos\alpha\n_1\otimes\n_1+\kappa_2\sin\alpha\n_1\otimes\n_2+\kappa_1\sin\alpha\n_2\otimes\n_1-\kappa_2\cos\alpha\n_2\otimes\n_2\}.
\end{equation}
To compute $\nablas\F$, we also need to express both $\nablas\n_1$
and $\nablas\n_2$ in the frame $\framen$,
\begin{equation}
	\label{eq:n_surface_gradients}
	\begin{cases}
		\nablas\n_1=\n_2\otimes\con-\kappa_1\normal\otimes\n_1,\\
		\nablas\n_2=-\n_1\otimes\con-\kappa_2\normal\otimes\n_2,
	\end{cases}	
\end{equation}
where $\con$ is the tangential field defined by \eqref{eq:connector_definition}.

A lengthy but simple calculation shows that equation \eqref{eq:integrability_tensor} acquires the following structure
\begin{equation}
	\label{eq:integrability_structure}
	\n_1\otimes\W_1+\n_2\otimes\W_2=\zero,
\end{equation}
where both $\W_1$ and $\W_2$ are skew-symmetric second-rank tensors, which by \eqref{eq:F_curvature_frame_representation}, \eqref{eq:n_surface_gradients}, and \eqref{eq:curvature_tensor_representation} are associated with the following axial vectors
\begin{equation}
	\label{eq:w_vectors}
	\begin{cases}
		\bm{w}_1=\nablas(\lambda\kappa_1\cos\alpha)\times\n_1+\nablas(\lambda\kappa_2\sin\alpha)\times\n_2-\lambda(\kappa_1+\kappa_2)\Rnu(\alpha)\n_2\times\con,\\
		\bm{w}_2=\nablas(\lambda\kappa_1\sin\alpha)\times\n_1-\nablas(\lambda\kappa_2\cos\alpha)\times\n_2-\lambda(\kappa_1+\kappa_2)\Rnu(\alpha)\n_1\times\con,
	\end{cases}
\end{equation}
which are both parallel to $\normal$. The system of scalar equations \eqref{eq:integrability_system} in the main text follows from requiring that $\bm{w}_1=\bm{w}_2\equiv\zero$ on $\surface$.

\begin{acknowledgements}
This paper was completed while AMS was visiting the Department of Mathematics of the University of Pavia, participating the University's Honours Programme ``Collegiale Non Residente'', where he taught a Master's course on \emph{Selected Topics in Fluid Dynamics}. The kind hospitality of the Department is gratefully acknowledged. EGV is a member of \emph{GNFM}, which is part of \emph{INdAM}, the Italian National Institute for Advanced Mathematics. This work was partly supported by \emph{GNFM}. We are grateful to Prof. R. Rosso for an educational discussion on the history of Scherk's surfaces. We also thank Prof. F. Burstall for an enlightening correspondence.
\end{acknowledgements}


\begin{thebibliography}{42}%
\makeatletter
\providecommand \@ifxundefined [1]{%
 \@ifx{#1\undefined}
}%
\providecommand \@ifnum [1]{%
 \ifnum #1\expandafter \@firstoftwo
 \else \expandafter \@secondoftwo
 \fi
}%
\providecommand \@ifx [1]{%
 \ifx #1\expandafter \@firstoftwo
 \else \expandafter \@secondoftwo
 \fi
}%
\providecommand \natexlab [1]{#1}%
\providecommand \enquote  [1]{``#1''}%
\providecommand \bibnamefont  [1]{#1}%
\providecommand \bibfnamefont [1]{#1}%
\providecommand \citenamefont [1]{#1}%
\providecommand \href@noop [0]{\@secondoftwo}%
\providecommand \href [0]{\begingroup \@sanitize@url \@href}%
\providecommand \@href[1]{\@@startlink{#1}\@@href}%
\providecommand \@@href[1]{\endgroup#1\@@endlink}%
\providecommand \@sanitize@url [0]{\catcode `\\12\catcode `\$12\catcode
  `\&12\catcode `\#12\catcode `\^12\catcode `\_12\catcode `\%12\relax}%
\providecommand \@@startlink[1]{}%
\providecommand \@@endlink[0]{}%
\providecommand \url  [0]{\begingroup\@sanitize@url \@url }%
\providecommand \@url [1]{\endgroup\@href {#1}{\urlprefix }}%
\providecommand \urlprefix  [0]{URL }%
\providecommand \Eprint [0]{\href }%
\providecommand \doibase [0]{https://doi.org/}%
\providecommand \selectlanguage [0]{\@gobble}%
\providecommand \bibinfo  [0]{\@secondoftwo}%
\providecommand \bibfield  [0]{\@secondoftwo}%
\providecommand \translation [1]{[#1]}%
\providecommand \BibitemOpen [0]{}%
\providecommand \bibitemStop [0]{}%
\providecommand \bibitemNoStop [0]{.\EOS\space}%
\providecommand \EOS [0]{\spacefactor3000\relax}%
\providecommand \BibitemShut  [1]{\csname bibitem#1\endcsname}%
\let\auto@bib@innerbib\@empty
\bibitem [{\citenamefont {Hildebrandt}\ and\ \citenamefont
  {Tromba}(1996)}]{hildebrandt:parsimonious}%
  \BibitemOpen
  \bibfield  {author} {\bibinfo {author} {\bibfnamefont {S.}~\bibnamefont
  {Hildebrandt}}\ and\ \bibinfo {author} {\bibfnamefont {A.}~\bibnamefont
  {Tromba}},\ }\href@noop {} {\emph {\bibinfo {title} {The Parsimonious
  Universe. {S}hape and Form in the Natural World}}}\ (\bibinfo  {publisher}
  {Springer-Verlag},\ \bibinfo {address} {New York},\ \bibinfo {year}
  {1996})\BibitemShut {NoStop}%
\bibitem [{\citenamefont {Antman}(1995)}]{antman:nonlinear}%
  \BibitemOpen
  \bibfield  {author} {\bibinfo {author} {\bibfnamefont {S.~S.}\ \bibnamefont
  {Antman}},\ }\href@noop {} {\emph {\bibinfo {title} {Nonlinear Problems of
  Elasticity}}},\ \bibinfo {series} {Applied Mathematical Sciences}, Vol.\
  \bibinfo {volume} {107}\ (\bibinfo  {publisher} {Springer},\ \bibinfo
  {address} {New York},\ \bibinfo {year} {1995})\BibitemShut {NoStop}%
\bibitem [{\citenamefont {Ghiba}\ \emph {et~al.}(2020)\citenamefont {Ghiba},
  \citenamefont {B\^irsan}, \citenamefont {Lewintan},\ and\ \citenamefont
  {Neff}}]{ghiba:isotropic}%
  \BibitemOpen
  \bibfield  {author} {\bibinfo {author} {\bibfnamefont {I.-D.}\ \bibnamefont
  {Ghiba}}, \bibinfo {author} {\bibfnamefont {M.}~\bibnamefont {B\^irsan}},
  \bibinfo {author} {\bibfnamefont {P.}~\bibnamefont {Lewintan}},\ and\
  \bibinfo {author} {\bibfnamefont {P.}~\bibnamefont {Neff}},\ }\bibfield
  {title} {\bibinfo {title} {The isotropic {C}osserat shell model including
  terms up to ${O}(h^5)$. {P}art {I}: {D}erivation in matrix notation},\ }\href
  {https://doi.org/https://doi.org/10.1007/s10659-020-09796-3} {\bibfield
  {journal} {\bibinfo  {journal} {J. Elast.}\ }\textbf {\bibinfo {volume}
  {142}},\ \bibinfo {pages} {201} (\bibinfo {year} {2020})}\BibitemShut
  {NoStop}%
\bibitem [{\citenamefont {Ghiba}\ \emph {et~al.}(2021)\citenamefont {Ghiba},
  \citenamefont {B\^irsan}, \citenamefont {Lewintan},\ and\ \citenamefont
  {Neff}}]{ghiba:constrained}%
  \BibitemOpen
  \bibfield  {author} {\bibinfo {author} {\bibfnamefont {I.-D.}\ \bibnamefont
  {Ghiba}}, \bibinfo {author} {\bibfnamefont {M.}~\bibnamefont {B\^irsan}},
  \bibinfo {author} {\bibfnamefont {P.}~\bibnamefont {Lewintan}},\ and\
  \bibinfo {author} {\bibfnamefont {P.}~\bibnamefont {Neff}},\ }\bibfield
  {title} {\bibinfo {title} {A constrained {C}osserat shell model up to order
  ${O}(h^5)$: {M}odelling, existence of minimizers, relations to classical
  shell models and scaling invariance of the bending tensor},\ }\href
  {https://doi.org/https://doi.org/10.1007/s10659-021-09851-7} {\bibfield
  {journal} {\bibinfo  {journal} {J. Elast.}\ }\textbf {\bibinfo {volume}
  {76}},\ \bibinfo {pages} {83} (\bibinfo {year} {2021})}\BibitemShut {NoStop}%
\bibitem [{\citenamefont {Vitral}\ and\ \citenamefont
  {Hanna}(2023{\natexlab{a}})}]{vitral:dilation}%
  \BibitemOpen
  \bibfield  {author} {\bibinfo {author} {\bibfnamefont {E.}~\bibnamefont
  {Vitral}}\ and\ \bibinfo {author} {\bibfnamefont {J.~A.}\ \bibnamefont
  {Hanna}},\ }\bibfield  {title} {\bibinfo {title} {Dilation-invariant bending
  of elastic plates, and broken symmetry in shells},\ }\href
  {https://doi.org/https://doi.org/10.1007/s10659-022-09894-4} {\bibfield
  {journal} {\bibinfo  {journal} {J. Elast.}\ }\textbf {\bibinfo {volume}
  {153}},\ \bibinfo {pages} {571} (\bibinfo {year}
  {2023}{\natexlab{a}})}\BibitemShut {NoStop}%
\bibitem [{\citenamefont {Vitral}\ and\ \citenamefont
  {Hanna}(2023{\natexlab{b}})}]{vitral:energies}%
  \BibitemOpen
  \bibfield  {author} {\bibinfo {author} {\bibfnamefont {E.}~\bibnamefont
  {Vitral}}\ and\ \bibinfo {author} {\bibfnamefont {J.~A.}\ \bibnamefont
  {Hanna}},\ }\bibfield  {title} {\bibinfo {title} {Energies for elastic plates
  and shells from quadratic-stretch elasticity},\ }\href
  {https://doi.org/https://doi.org/10.1007/s10659-022-09895-3} {\bibfield
  {journal} {\bibinfo  {journal} {J. Elast.}\ }\textbf {\bibinfo {volume}
  {153}},\ \bibinfo {pages} {581} (\bibinfo {year}
  {2023}{\natexlab{b}})}\BibitemShut {NoStop}%
\bibitem [{\citenamefont {Virga}(2024)}]{virga:pure}%
  \BibitemOpen
  \bibfield  {author} {\bibinfo {author} {\bibfnamefont {E.~G.}\ \bibnamefont
  {Virga}},\ }\bibfield  {title} {\bibinfo {title} {Pure measures of bending
  for soft plates},\ }\href {https://doi.org/10.1039/D3SM01123B} {\bibfield
  {journal} {\bibinfo  {journal} {Soft Matter}\ }\textbf {\bibinfo {volume}
  {20}},\ \bibinfo {pages} {144} (\bibinfo {year} {2024})}\BibitemShut
  {NoStop}%
\bibitem [{\citenamefont {Acharya}(2024)}]{acharya:mid-surface}%
  \BibitemOpen
  \bibfield  {author} {\bibinfo {author} {\bibfnamefont {A.}~\bibnamefont
  {Acharya}},\ }\bibfield  {title} {\bibinfo {title} {Mid-surface scaling
  invariance of some bending strain measures},\ }\href
  {https://doi.org/https://doi.org/10.1007/s10659-024-10066-9} {\bibfield
  {journal} {\bibinfo  {journal} {J. Elast.}\ ,\ \bibinfo {pages} {5517}}
  (\bibinfo {year} {2024})}\BibitemShut {NoStop}%
\bibitem [{\citenamefont {Ghiba}\ \emph {et~al.}(2023)\citenamefont {Ghiba},
  \citenamefont {Lewintan}, \citenamefont {Sky},\ and\ \citenamefont
  {Neff}}]{ghiba:essay}%
  \BibitemOpen
  \bibfield  {author} {\bibinfo {author} {\bibfnamefont {I.-D.}\ \bibnamefont
  {Ghiba}}, \bibinfo {author} {\bibfnamefont {P.}~\bibnamefont {Lewintan}},
  \bibinfo {author} {\bibfnamefont {A.}~\bibnamefont {Sky}},\ and\ \bibinfo
  {author} {\bibfnamefont {P.}~\bibnamefont {Neff}},\ }\href
  {https://doi.org/https://doi.org/10.48550/arXiv.2312.10928} {\bibinfo {title}
  {An essay on deformation measures in isotropic thin shell theories. bending
  versus curvature}} (\bibinfo {year} {2023}),\ \Eprint
  {https://arxiv.org/abs/2312.10928} {arXiv:2312.10928 [math-ph]} \BibitemShut
  {NoStop}%
\bibitem [{\citenamefont {Vitral}\ and\ \citenamefont
  {Hanna}(2024)}]{vitral:assorted}%
  \BibitemOpen
  \bibfield  {author} {\bibinfo {author} {\bibfnamefont {E.}~\bibnamefont
  {Vitral}}\ and\ \bibinfo {author} {\bibfnamefont {J.~A.}\ \bibnamefont
  {Hanna}},\ }\href {https://doi.org/https://doi.org/10.48550/arXiv.2405.06638}
  {\bibinfo {title} {Assorted remarks on bending measures and energies for
  plates and shells, and their invariance properties}} (\bibinfo {year}
  {2024}),\ \Eprint {https://arxiv.org/abs/2405.06638} {arXiv:2405.06638
  [cond-mat.soft]} \BibitemShut {NoStop}%
\bibitem [{\citenamefont
  {Weatherburn}(2016{\natexlab{a}})}]{weatherburn:differential_1}%
  \BibitemOpen
  \bibfield  {author} {\bibinfo {author} {\bibfnamefont {C.~E.}\ \bibnamefont
  {Weatherburn}},\ }\href@noop {} {\emph {\bibinfo {title} {Differential
  Geometry of Three Dimensions}}},\ Vol.~\bibinfo {volume} {I}\ (\bibinfo
  {publisher} {Cambridge University Press},\ \bibinfo {address} {Cambridge},\
  \bibinfo {year} {2016})\BibitemShut {NoStop}%
\bibitem [{\citenamefont {\v{S}ilhav\'y}(2021)}]{silhavy:new}%
  \BibitemOpen
  \bibfield  {author} {\bibinfo {author} {\bibfnamefont {M.}~\bibnamefont
  {\v{S}ilhav\'y}},\ }\bibfield  {title} {\bibinfo {title} {A new approach to
  curvature measures in linear shell theories},\ }\href
  {https://doi.org/10.1177/1081286520972752} {\bibfield  {journal} {\bibinfo
  {journal} {Math. Mech. Solids}\ }\textbf {\bibinfo {volume} {26}},\ \bibinfo
  {pages} {1241} (\bibinfo {year} {2021})}\BibitemShut {NoStop}%
\bibitem [{\citenamefont {Kralj}\ \emph {et~al.}(2011)\citenamefont {Kralj},
  \citenamefont {Rosso},\ and\ \citenamefont {Virga}}]{kralj:curvature}%
  \BibitemOpen
  \bibfield  {author} {\bibinfo {author} {\bibfnamefont {S.}~\bibnamefont
  {Kralj}}, \bibinfo {author} {\bibfnamefont {R.}~\bibnamefont {Rosso}},\ and\
  \bibinfo {author} {\bibfnamefont {E.~G.}\ \bibnamefont {Virga}},\ }\bibfield
  {title} {\bibinfo {title} {Curvature control of valence on nematic shells},\
  }\href {https://doi.org/10.1039/C0SM00378F} {\bibfield  {journal} {\bibinfo
  {journal} {Soft Matter}\ }\textbf {\bibinfo {volume} {7}},\ \bibinfo {pages}
  {670} (\bibinfo {year} {2011})}\BibitemShut {NoStop}%
\bibitem [{\citenamefont {Rosso}\ \emph {et~al.}(2012)\citenamefont {Rosso},
  \citenamefont {Virga},\ and\ \citenamefont {Kralj}}]{rosso:parallel}%
  \BibitemOpen
  \bibfield  {author} {\bibinfo {author} {\bibfnamefont {R.}~\bibnamefont
  {Rosso}}, \bibinfo {author} {\bibfnamefont {E.~G.}\ \bibnamefont {Virga}},\
  and\ \bibinfo {author} {\bibfnamefont {S.}~\bibnamefont {Kralj}},\ }\bibfield
   {title} {\bibinfo {title} {Parallel transport and defects on nematic
  shells},\ }\href {https://doi.org/https://doi.org/10.1007/s00161-012-0259-4}
  {\bibfield  {journal} {\bibinfo  {journal} {Continuum Mech. Thermodyn.}\
  }\textbf {\bibinfo {volume} {24}},\ \bibinfo {pages} {643} (\bibinfo {year}
  {2012})}\BibitemShut {NoStop}%
\bibitem [{\citenamefont
  {Weatherburn}(2016{\natexlab{b}})}]{weatherburn:differential_2}%
  \BibitemOpen
  \bibfield  {author} {\bibinfo {author} {\bibfnamefont {C.~E.}\ \bibnamefont
  {Weatherburn}},\ }\href@noop {} {\emph {\bibinfo {title} {Differential
  Geometry of Three Dimensions}}},\ Vol.~\bibinfo {volume} {II}\ (\bibinfo
  {publisher} {Cambridge University Press},\ \bibinfo {address} {Cambridge},\
  \bibinfo {year} {2016})\BibitemShut {NoStop}%
\bibitem [{\citenamefont {Beltrami}(1902)}]{beltrami:opere}%
  \BibitemOpen
  \bibfield  {author} {\bibinfo {author} {\bibfnamefont {E.}~\bibnamefont
  {Beltrami}},\ }\href@noop {} {\emph {\bibinfo {title} {Opere Matematiche}}},\
  Vol.~\bibinfo {volume} {1}\ (\bibinfo  {publisher} {Hoepli},\ \bibinfo
  {address} {Milan},\ \bibinfo {year} {1902})\ \bibinfo {note} {available from
  \url{https://gallica.bnf.fr/ark:/12148/bpt6k99432q/f1.item}}\BibitemShut
  {NoStop}%
\bibitem [{\citenamefont {Gurtin}\ \emph {et~al.}(2010)\citenamefont {Gurtin},
  \citenamefont {Fried},\ and\ \citenamefont {Anand}}]{gurtin:mechanics}%
  \BibitemOpen
  \bibfield  {author} {\bibinfo {author} {\bibfnamefont {M.~E.}\ \bibnamefont
  {Gurtin}}, \bibinfo {author} {\bibfnamefont {E.}~\bibnamefont {Fried}},\ and\
  \bibinfo {author} {\bibfnamefont {L.}~\bibnamefont {Anand}},\ }\href@noop {}
  {\emph {\bibinfo {title} {The Mechanics and Thermodynamics of Continua}}}\
  (\bibinfo  {publisher} {Cambridge University Press},\ \bibinfo {address}
  {Cambridge},\ \bibinfo {year} {2010})\BibitemShut {NoStop}%
\bibitem [{\citenamefont {Man}\ and\ \citenamefont
  {Cohen}(1986)}]{man:coordinate}%
  \BibitemOpen
  \bibfield  {author} {\bibinfo {author} {\bibfnamefont {C.-S.}\ \bibnamefont
  {Man}}\ and\ \bibinfo {author} {\bibfnamefont {H.}~\bibnamefont {Cohen}},\
  }\bibfield  {title} {\bibinfo {title} {A coordinate-free approach to the
  kinematics of membranes},\ }\href
  {https://doi.org/https://doi.org/10.1007/BF00041068} {\bibfield  {journal}
  {\bibinfo  {journal} {J. Elast.}\ }\textbf {\bibinfo {volume} {16}},\
  \bibinfo {pages} {97} (\bibinfo {year} {1986})}\BibitemShut {NoStop}%
\bibitem [{\citenamefont {Gurtin}\ and\ \citenamefont
  {Murdoch}(1975{\natexlab{a}})}]{gurtin:continuum}%
  \BibitemOpen
  \bibfield  {author} {\bibinfo {author} {\bibfnamefont {M.~E.}\ \bibnamefont
  {Gurtin}}\ and\ \bibinfo {author} {\bibfnamefont {A.~I.}\ \bibnamefont
  {Murdoch}},\ }\bibfield  {title} {\bibinfo {title} {A continuum theory of
  elastic material surfaces},\ }\href
  {https://doi.org/https://doi.org/10.1007/BF00261375} {\bibfield  {journal}
  {\bibinfo  {journal} {Arch. Rational Mech. Anal.}\ }\textbf {\bibinfo
  {volume} {57}},\ \bibinfo {pages} {291} (\bibinfo {year}
  {1975}{\natexlab{a}})},\ \bibinfo {note} {see also
  \cite{gurtin:addenda}}\BibitemShut {NoStop}%
\bibitem [{\citenamefont {Pietraszkiewicz}\ \emph {et~al.}(2008)\citenamefont
  {Pietraszkiewicz}, \citenamefont {Szwabowicz},\ and\ \citenamefont
  {Vall\'ee}}]{pietraszkiewicz:determination}%
  \BibitemOpen
  \bibfield  {author} {\bibinfo {author} {\bibfnamefont {W.}~\bibnamefont
  {Pietraszkiewicz}}, \bibinfo {author} {\bibfnamefont {M.}~\bibnamefont
  {Szwabowicz}},\ and\ \bibinfo {author} {\bibfnamefont {C.}~\bibnamefont
  {Vall\'ee}},\ }\bibfield  {title} {\bibinfo {title} {Determination of the
  midsurface of a deformed shell from prescribed surface strains and bendings
  via the polar decomposition},\ }\href
  {https://doi.org/https://doi.org/10.1016/j.ijnonlinmec.2008.02.003}
  {\bibfield  {journal} {\bibinfo  {journal} {Int. J. Non-Linear Mech.}\
  }\textbf {\bibinfo {volume} {43}},\ \bibinfo {pages} {579} (\bibinfo {year}
  {2008})}\BibitemShut {NoStop}%
\bibitem [{\citenamefont {Palais}\ \emph {et~al.}(2009)\citenamefont {Palais},
  \citenamefont {Palais},\ and\ \citenamefont {Rodi}}]{palais:disorienting}%
  \BibitemOpen
  \bibfield  {author} {\bibinfo {author} {\bibfnamefont {B.}~\bibnamefont
  {Palais}}, \bibinfo {author} {\bibfnamefont {R.}~\bibnamefont {Palais}},\
  and\ \bibinfo {author} {\bibfnamefont {S.}~\bibnamefont {Rodi}},\ }\bibfield
  {title} {\bibinfo {title} {A disorienting look at {E}uler's theorem on the
  axis of a rotation},\ }\href {https://doi.org/10.4169/000298909X477014}
  {\bibfield  {journal} {\bibinfo  {journal} {Am. Math. Mon.}\ }\textbf
  {\bibinfo {volume} {116}},\ \bibinfo {pages} {892} (\bibinfo {year}
  {2009})}\BibitemShut {NoStop}%
\bibitem [{\citenamefont {Hughes}\ and\ \citenamefont
  {Brezzi}(1989)}]{hughes:drilling}%
  \BibitemOpen
  \bibfield  {author} {\bibinfo {author} {\bibfnamefont {T.~J.}\ \bibnamefont
  {Hughes}}\ and\ \bibinfo {author} {\bibfnamefont {F.}~\bibnamefont
  {Brezzi}},\ }\bibfield  {title} {\bibinfo {title} {On drilling degrees of
  freedom},\ }\href
  {https://doi.org/https://doi.org/10.1016/0045-7825(89)90124-2} {\bibfield
  {journal} {\bibinfo  {journal} {Comput. Methods Appl. Mech. Engrg.}\ }\textbf
  {\bibinfo {volume} {72}},\ \bibinfo {pages} {105} (\bibinfo {year}
  {1989})}\BibitemShut {NoStop}%
\bibitem [{\citenamefont {Fox}\ and\ \citenamefont {Simo}(1992)}]{fox:drill}%
  \BibitemOpen
  \bibfield  {author} {\bibinfo {author} {\bibfnamefont {D.}~\bibnamefont
  {Fox}}\ and\ \bibinfo {author} {\bibfnamefont {J.}~\bibnamefont {Simo}},\
  }\bibfield  {title} {\bibinfo {title} {A drill rotation formulation for
  geometrically exact shells},\ }\href
  {https://doi.org/https://doi.org/10.1016/0045-7825(92)90002-2} {\bibfield
  {journal} {\bibinfo  {journal} {Comput. Methods Appl. Mech. Engrg.}\ }\textbf
  {\bibinfo {volume} {98}},\ \bibinfo {pages} {329} (\bibinfo {year}
  {1992})}\BibitemShut {NoStop}%
\bibitem [{\citenamefont {Mohammadi~Saem}\ \emph {et~al.}(2021)\citenamefont
  {Mohammadi~Saem}, \citenamefont {Lewintan},\ and\ \citenamefont
  {Neff}}]{saem:in-plane}%
  \BibitemOpen
  \bibfield  {author} {\bibinfo {author} {\bibfnamefont {M.}~\bibnamefont
  {Mohammadi~Saem}}, \bibinfo {author} {\bibfnamefont {P.}~\bibnamefont
  {Lewintan}},\ and\ \bibinfo {author} {\bibfnamefont {P.}~\bibnamefont
  {Neff}},\ }\bibfield  {title} {\bibinfo {title} {On in-plane drill rotations
  for cosserat surfaces},\ }\href {https://doi.org/10.1098/rspa.2021.0158}
  {\bibfield  {journal} {\bibinfo  {journal} {Proc. R. Soc. Lond. A}\ }\textbf
  {\bibinfo {volume} {477}},\ \bibinfo {pages} {20210158} (\bibinfo {year}
  {2021})}\BibitemShut {NoStop}%
\bibitem [{\citenamefont {Altmann}(1989)}]{altmann:hamilton}%
  \BibitemOpen
  \bibfield  {author} {\bibinfo {author} {\bibfnamefont {S.~L.}\ \bibnamefont
  {Altmann}},\ }\bibfield  {title} {\bibinfo {title} {Hamilton, {R}odrigues,
  and the quaternion scandal},\ }\href
  {https://doi.org/10.1080/0025570X.1989.11977459} {\bibfield  {journal}
  {\bibinfo  {journal} {Math. Mag.}\ }\textbf {\bibinfo {volume} {62}},\
  \bibinfo {pages} {291} (\bibinfo {year} {1989})}\BibitemShut {NoStop}%
\bibitem [{\citenamefont {Rodrigues}(1840)}]{rodrigues:lois}%
  \BibitemOpen
  \bibfield  {author} {\bibinfo {author} {\bibfnamefont {O.}~\bibnamefont
  {Rodrigues}},\ }\bibfield  {title} {\bibinfo {title} {Des lois
  g\'eom\'etriques qui r\'egissent les d\'eplacements d'un syst\`eme solide
  dans l'espace, et de la variation des coordonn\'ees provenant de ces
  d\'eplacements consid\'er\'es ind\'ependamment des causes qui peuvent les
  produire},\ }\href@noop {} {\bibfield  {journal} {\bibinfo  {journal} {J.
  Math. Pures Appl.}\ }\textbf {\bibinfo {volume} {5}},\ \bibinfo {pages} {380}
  (\bibinfo {year} {1840})},\ \bibinfo {note} {available from
  \url{http://www.numdam.org/item/JMPA_1840_1_5__380_0/}}\BibitemShut {NoStop}%
\bibitem [{\citenamefont {Mladenova}\ and\ \citenamefont
  {Mladenov}(2011)}]{mladenova:vector}%
  \BibitemOpen
  \bibfield  {author} {\bibinfo {author} {\bibfnamefont {C.~D.}\ \bibnamefont
  {Mladenova}}\ and\ \bibinfo {author} {\bibfnamefont {I.~M.}\ \bibnamefont
  {Mladenov}},\ }\bibfield  {title} {\bibinfo {title} {Vector decomposition of
  finite rotations},\ }\href
  {https://doi.org/https://doi.org/10.1016/S0034-4877(11)60030-X} {\bibfield
  {journal} {\bibinfo  {journal} {Rep. Math. Phys.}\ }\textbf {\bibinfo
  {volume} {68}},\ \bibinfo {pages} {107} (\bibinfo {year} {2011})}\BibitemShut
  {NoStop}%
\bibitem [{\citenamefont {Martins}\ and\ \citenamefont
  {Podio-Guidugli}(1979)}]{martins:variational}%
  \BibitemOpen
  \bibfield  {author} {\bibinfo {author} {\bibfnamefont {L.~C.}\ \bibnamefont
  {Martins}}\ and\ \bibinfo {author} {\bibfnamefont {P.}~\bibnamefont
  {Podio-Guidugli}},\ }\bibfield  {title} {\bibinfo {title} {A variational
  approach to the polar decomposition theorem},\ }\href@noop {} {\bibfield
  {journal} {\bibinfo  {journal} {Atti Accad. Naz. Lincei Cl. Sci. Fis. Mat.
  Natur. Rend.}\ }\textbf {\bibinfo {volume} {66}},\ \bibinfo {pages} {487}
  (\bibinfo {year} {1979})},\ \bibinfo {note} {available from
  \url{http://www.bdim.eu/item?id=RLINA_1979_8_66_6_487_0}}\BibitemShut
  {NoStop}%
\bibitem [{\citenamefont {Martins}\ and\ \citenamefont
  {Podio-Guidugli}(1980)}]{martins:variational_1980}%
  \BibitemOpen
  \bibfield  {author} {\bibinfo {author} {\bibfnamefont {L.~C.}\ \bibnamefont
  {Martins}}\ and\ \bibinfo {author} {\bibfnamefont {P.}~\bibnamefont
  {Podio-Guidugli}},\ }\bibfield  {title} {\bibinfo {title} {An elementary
  proof of the polar decomposition theorem},\ }\href
  {https://doi.org/10.1080/00029890.1980.11995017} {\bibfield  {journal}
  {\bibinfo  {journal} {Amer. Math. Month.}\ }\textbf {\bibinfo {volume}
  {87}},\ \bibinfo {pages} {288} (\bibinfo {year} {1980})}\BibitemShut
  {NoStop}%
\bibitem [{\citenamefont {Grioli}(1940)}]{grioli:proprieta}%
  \BibitemOpen
  \bibfield  {author} {\bibinfo {author} {\bibfnamefont {G.}~\bibnamefont
  {Grioli}},\ }\bibfield  {title} {\bibinfo {title} {Una propriet\`a di minimo
  nella cinematica delle deformazioni finite},\ }\href@noop {} {\bibfield
  {journal} {\bibinfo  {journal} {Boll. Un. Math. Ital.}\ }\textbf {\bibinfo
  {volume} {2}},\ \bibinfo {pages} {252} (\bibinfo {year} {1940})},\ \bibinfo
  {note} {available from
  \url{https://www.uni-due.de//imperia/md/content/mathematik/ag_neff/grioli_deformatione.pdf}}\BibitemShut
  {NoStop}%
\bibitem [{\citenamefont {Neff}\ \emph {et~al.}(2014)\citenamefont {Neff},
  \citenamefont {Lankeit},\ and\ \citenamefont {Madeo}}]{neff:grioli's}%
  \BibitemOpen
  \bibfield  {author} {\bibinfo {author} {\bibfnamefont {P.}~\bibnamefont
  {Neff}}, \bibinfo {author} {\bibfnamefont {J.}~\bibnamefont {Lankeit}},\ and\
  \bibinfo {author} {\bibfnamefont {A.}~\bibnamefont {Madeo}},\ }\bibfield
  {title} {\bibinfo {title} {On {G}rioli’s minimum property and its relation
  to {C}auchy’s polar decomposition},\ }\href
  {https://doi.org/https://doi.org/10.1016/j.ijengsci.2014.02.026} {\bibfield
  {journal} {\bibinfo  {journal} {Int. J. Eng. Sci.}\ }\textbf {\bibinfo
  {volume} {80}},\ \bibinfo {pages} {209} (\bibinfo {year} {2014})}\BibitemShut
  {NoStop}%
\bibitem [{\citenamefont {Ozenda}\ \emph {et~al.}(2020)\citenamefont {Ozenda},
  \citenamefont {Sonnet},\ and\ \citenamefont {Virga}}]{ozenda:blend}%
  \BibitemOpen
  \bibfield  {author} {\bibinfo {author} {\bibfnamefont {O.}~\bibnamefont
  {Ozenda}}, \bibinfo {author} {\bibfnamefont {A.~M.}\ \bibnamefont {Sonnet}},\
  and\ \bibinfo {author} {\bibfnamefont {E.~G.}\ \bibnamefont {Virga}},\
  }\bibfield  {title} {\bibinfo {title} {A blend of stretching and bending in
  nematic polymer networks},\ }\href
  {https://doi.org/https://doi.org/10.1039/D0SM00642D} {\bibfield  {journal}
  {\bibinfo  {journal} {Soft Matter}\ }\textbf {\bibinfo {volume} {16}},\
  \bibinfo {pages} {8877} (\bibinfo {year} {2020})}\BibitemShut {NoStop}%
\bibitem [{\citenamefont {Ozenda}\ and\ \citenamefont
  {Virga}(2021)}]{ozenda:kirchhoff}%
  \BibitemOpen
  \bibfield  {author} {\bibinfo {author} {\bibfnamefont {O.}~\bibnamefont
  {Ozenda}}\ and\ \bibinfo {author} {\bibfnamefont {E.~G.}\ \bibnamefont
  {Virga}},\ }\bibfield  {title} {\bibinfo {title} {On the {K}irchhoff-{L}ove
  hypothesis (revised and vindicated)},\ }\href
  {https://doi.org/https://doi.org/10.1007/s10659-021-09819-7} {\bibfield
  {journal} {\bibinfo  {journal} {J. Elast.}\ }\textbf {\bibinfo {volume}
  {143}},\ \bibinfo {pages} {359} (\bibinfo {year} {2021})}\BibitemShut
  {NoStop}%
\bibitem [{\citenamefont {Dierkes}\ \emph {et~al.}(2010)\citenamefont
  {Dierkes}, \citenamefont {Hildebrandt},\ and\ \citenamefont
  {Sauvigny}}]{dierkes:minimal}%
  \BibitemOpen
  \bibfield  {author} {\bibinfo {author} {\bibfnamefont {U.}~\bibnamefont
  {Dierkes}}, \bibinfo {author} {\bibfnamefont {S.}~\bibnamefont
  {Hildebrandt}},\ and\ \bibinfo {author} {\bibfnamefont {F.}~\bibnamefont
  {Sauvigny}},\ }\href@noop {} {\emph {\bibinfo {title} {Minimal Surfaces}}},\
  \bibinfo {edition} {2nd}\ ed.,\ \bibinfo {series} {A Series of Comprehensive
  Studies in Mathematics}, Vol.\ \bibinfo {volume} {339}\ (\bibinfo
  {publisher} {Springer},\ \bibinfo {address} {Berlin},\ \bibinfo {year}
  {2010})\BibitemShut {NoStop}%
\bibitem [{\citenamefont {Chern}(1955)}]{chern:elementary}%
  \BibitemOpen
  \bibfield  {author} {\bibinfo {author} {\bibfnamefont {S.-S.}\ \bibnamefont
  {Chern}},\ }\bibfield  {title} {\bibinfo {title} {An elementary proof of the
  existence of isothermal parameters on a surface},\ }\href
  {https://doi.org/https://doi.org/10.1090/S0002-9939-1955-0074856-1}
  {\bibfield  {journal} {\bibinfo  {journal} {Proc. Amer. Math. Soc.}\ }\textbf
  {\bibinfo {volume} {6}},\ \bibinfo {pages} {771} (\bibinfo {year}
  {1955})}\BibitemShut {NoStop}%
\bibitem [{\citenamefont {Osserman}(1986)}]{osserman:survey}%
  \BibitemOpen
  \bibfield  {author} {\bibinfo {author} {\bibfnamefont {R.}~\bibnamefont
  {Osserman}},\ }\href@noop {} {\emph {\bibinfo {title} {A Survey of Minimal
  Surfaces}}}\ (\bibinfo  {publisher} {Dover},\ \bibinfo {address} {Mineola,
  NY},\ \bibinfo {year} {1986})\BibitemShut {NoStop}%
\bibitem [{\citenamefont {Spivak}(1999)}]{spivack:comprehensive_4}%
  \BibitemOpen
  \bibfield  {author} {\bibinfo {author} {\bibfnamefont {M.}~\bibnamefont
  {Spivak}},\ }\href@noop {} {\emph {\bibinfo {title} {A Comprehensive
  Introduction to Differential Geometry}}},\ \bibinfo {edition} {3rd}\ ed.,\
  Vol.~\bibinfo {volume} {4}\ (\bibinfo  {publisher} {Publish or Perish},\
  \bibinfo {address} {Houston,Texas},\ \bibinfo {year} {1999})\BibitemShut
  {NoStop}%
\bibitem [{\citenamefont
  {Weierstrass}(1903)}]{weierstrass:untersuchungen_berlin}%
  \BibitemOpen
  \bibfield  {author} {\bibinfo {author} {\bibfnamefont {K.}~\bibnamefont
  {Weierstrass}},\ }\bibfield  {title} {\bibinfo {title} {Untersuchungen \"uber
  die {F}l\"achen, deren mittlere {K}r\"ummung \"uberall gleich {N}ull ist},\
  }in\ \href@noop {} {\emph {\bibinfo {booktitle} {Mathematische {W}erke:
  {H}erausgegeben unter {M}itwirkung einer von der k\"oniglich preussischen
  {A}kademie der {W}issenschaften eingesetzten {C}ommission}}},\ \bibinfo
  {editor} {edited by\ \bibinfo {editor} {\bibfnamefont {J.}~\bibnamefont
  {Knoblauch}}}\ (\bibinfo  {publisher} {Mayer \& M\"uller},\ \bibinfo
  {address} {Berlin},\ \bibinfo {year} {1903})\ pp.\ \bibinfo {pages}
  {39--52},\ \bibinfo {note} {available from
  \url{https://archive.org/details/mathematischewer03weieuoft/page/38/mode/2up}}\BibitemShut
  {NoStop}%
\bibitem [{\citenamefont
  {Weierstrass}(2013)}]{weierstrass:untersuchungen_cambridge}%
  \BibitemOpen
  \bibfield  {author} {\bibinfo {author} {\bibfnamefont {K.}~\bibnamefont
  {Weierstrass}},\ }\bibfield  {title} {\bibinfo {title} {Untersuchungen \"uber
  die {F}l\"achen, deren mittlere {K}r\"ummung \"uberall gleich {N}ull ist},\
  }in\ \href@noop {} {\emph {\bibinfo {booktitle} {Mathematische {W}erke:
  {H}erausgegeben unter {M}itwirkung einer von der k\"oniglich preussischen
  {A}kademie der {W}issenschaften eingesetzten {C}ommission}}},\ \bibinfo
  {series and number} {Cambridge Library Collection -- Mathematics},\ \bibinfo
  {editor} {edited by\ \bibinfo {editor} {\bibfnamefont {J.}~\bibnamefont
  {Knoblauch}}}\ (\bibinfo  {publisher} {Cambridge University Press},\ \bibinfo
  {address} {Cambridge},\ \bibinfo {year} {2013})\ pp.\ \bibinfo {pages}
  {39--52}\BibitemShut {NoStop}%
\bibitem [{\citenamefont {Bour}(1862)}]{bour:theorie}%
  \BibitemOpen
  \bibfield  {author} {\bibinfo {author} {\bibfnamefont {E.}~\bibnamefont
  {Bour}},\ }\bibfield  {title} {\bibinfo {title} {Th{\'e}orie de la
  d{\'e}formation des surfaces},\ }\href@noop {} {\bibfield  {journal}
  {\bibinfo  {journal} {J. {\'E}c. Imp{\'e}riale Polytech.}\ }\textbf {\bibinfo
  {volume} {39}},\ \bibinfo {pages} {1} (\bibinfo {year} {1862})}\BibitemShut
  {NoStop}%
\bibitem [{\citenamefont {Cogliati}\ and\ \citenamefont
  {Rivis}(2022)}]{cogliati:origins}%
  \BibitemOpen
  \bibfield  {author} {\bibinfo {author} {\bibfnamefont {A.}~\bibnamefont
  {Cogliati}}\ and\ \bibinfo {author} {\bibfnamefont {R.}~\bibnamefont
  {Rivis}},\ }\bibfield  {title} {\bibinfo {title} {The origins of the
  fundamental theorem of surface theory},\ }\href
  {https://doi.org/https://doi.org/10.1016/j.hm.2022.09.001} {\bibfield
  {journal} {\bibinfo  {journal} {Hist. Math.}\ }\textbf {\bibinfo {volume}
  {61}},\ \bibinfo {pages} {45} (\bibinfo {year} {2022})}\BibitemShut {NoStop}%
\bibitem [{\citenamefont {Gurtin}\ and\ \citenamefont
  {Murdoch}(1975{\natexlab{b}})}]{gurtin:addenda}%
  \BibitemOpen
  \bibfield  {author} {\bibinfo {author} {\bibfnamefont {M.~E.}\ \bibnamefont
  {Gurtin}}\ and\ \bibinfo {author} {\bibfnamefont {A.~I.}\ \bibnamefont
  {Murdoch}},\ }\bibfield  {title} {\bibinfo {title} {Addenda to our paper {A}
  continuum theory of elastic material surfaces},\ }\href
  {https://doi.org/https://doi.org/10.1007/BF00250426} {\bibfield  {journal}
  {\bibinfo  {journal} {Arch. Rational Mech. Anal.}\ }\textbf {\bibinfo
  {volume} {59}},\ \bibinfo {pages} {389} (\bibinfo {year}
  {1975}{\natexlab{b}})}\BibitemShut {NoStop}%
\end{thebibliography}

\begin{thebibliography}{9}

\bibitem{1} Allwood JM, Cullen JM. 2011 \textit{Sustainable materials:  with both eyes open}.
Cambridge, UK: UIT Cambridge. See \href{http://www.withbotheyesopen.com}{http://www.withbotheyesopen.com}.

\bibitem{2}  MacKay DJC. 2008  \textit{Sustainable energy:  without the hot air}.
 Cambridge, UK: UIT Cambridge. See \href{http://www.withouthotair.com}{http://www.withouthotair.com}.

\bibitem{3} Gallman PG. 2011  \textit{Green alternatives and national energy strategy: the facts
 behind the headlines}.  Baltimore,\ MD: Johns Hopkins University Press.

\bibitem{4} MacKay DJC. 2013.  Solar energy in the context of energy use, energy transportation, and
 energy storage. \textit{Proc. R. Soc. A} \textbf{371}.

\end{thebibliography}
\end{document}